\documentclass[12pt,a4paper]{amsart}
\usepackage[T2A]{fontenc}
\usepackage[cp1251]{inputenc}
\usepackage[english]{babel}
\usepackage{amssymb, amsmath, latexsym}
\usepackage{amsfonts,mathtext,cite,enumerate,float,amsthm}
\numberwithin{equation}{section}
\renewcommand{\le}{\leqslant}
\renewcommand{\ge}{\geqslant}
\newcommand{\rank}{\operatorname{rank}}

\newenvironment{Proof}
{{\bf Proof.}}
{\hfill$\scriptstyle\blacksquare$}
\newtheorem{theorem}{Theorem}
\newtheorem*{thL}{Theorem}
\newtheorem{lemma}{Lemma}
\theoremstyle{definition}
\newtheorem{remark}{Remark}
\newtheorem{example}{Example}

\usepackage{geometry}
\begin{document}
\title[Approximation by amplitude and frequency operators]{Approximation by\\ amplitude and frequency operators}
\author{Petr Chunaev and Vladimir Danchenko}
\footnotetext{The work of P. Chunaev was financially supported by RFBR (the projects No. 12-01-31471
mol\_a and No. 16-31-00252 mol\_a). The work of V. Danchenko was financially supported by
DRPNNiT No. 1.1348.2011, RFBR (the project No. 14-01-00510), the Russian Ministry of
Education (No. 2014/13, the project code 3037; No. 1.574.2016/FPM) and carried out within
the Vladimir State University state task No. 2014/13 in the field of scientific activity (the topic
2868).}

\maketitle

\begin{abstract}
We study Pad\'{e} interpolation at the node $z=0$ of functions $f(z)=\sum_{m=0}^{\infty} f_m z^m$, analytic in a neighbourhood of this node,  by {\it amplitude and frequency operators} (\textit{sums}) of the form
$$
 \sum_{k=1}^n \mu_k
h(\lambda_k z), \qquad \mu_k,\lambda_k\in \mathbb{C}.
$$
Here $h(z)=\sum_{m=0}^{\infty} h_m z^m$, $h_m\ne 0$, is a fixed (\textit{basis}) function, analytic at the origin, and the interpolation is carried out by an appropriate choice of \textit{amplitudes} $\mu_k $ and \textit{frequencies}
$\lambda_k$. The solvability of the $2n$-multiple interpolation problem is determined
by the solvability of the associated moment problem
$$
\sum_{k=1}^n\mu_k \lambda_k^m={f_m}/{h_m}, \qquad
m=\overline{0,2n-1}.
$$

In a number of cases, when the moment problem is consistent, it can be solved by the
classical method due to Prony and Sylvester, moreover, one can easily construct the corresponding interpolating sum too.
In the case of inconsistent moment problems, we propose a regularization method, which consists in adding a special
binomial $c_1z^{n-1}+c_2 z^{2n-1}$ to an amplitude and frequency
sum so that the moment problem, associated with the sum obtained, can be already solved by the method of Prony and Sylvester.
This approach enables us to obtain interpolation formulas with $n$
nodes $\lambda_k z$, being exact for the polynomials of degree $\le 2n-1$, whilst traditional
formulas with the same number of nodes are usually exact only for
the polynomials of degree $\le n-1$. The regularization method is applied to numerical differentiation and extrapolation.
\end{abstract}

\section{Introduction and statement of the problem}

In \cite{Dan2008,DanChu2011,Chu2010,Chu2012} the so-called  \textit{$h$-sums} of the form
\begin{equation*}
\label{h-sum}
\mathcal{H}_n(\{\lambda_k\},h;z)=\sum_{k=1}^n\lambda_k h(\lambda_k
z), \qquad z, \lambda_k\in \mathbb{C},\qquad n\in \mathbb{N},
\end{equation*}
are studied. Hereinafter $h(z)=\sum_{m=0}^\infty h_mz^m$ is a
function, analytic in a disc  $|z|<\rho$, $\rho>0$. We call it \textit{a basis function}. Obviously, $\mathcal{H}_n(z)$ is
well-defined and analytic in the disc $|z|<\rho\cdot
\min_{k=\overline{1,n}} |\lambda_k|^{-1}$.

In \cite{Dan2008,DanChu2011,Chu2010,Chu2012} operators
$\mathcal{H}_n(\{\lambda_k\},h;z)$ are used as a tool for
$n$-multiple (Pad\'{e}) interpolation and approximation of functions
$f$, analytic in a neighbourhood of the origin. In
particular, it is shown in \cite{Dan2008} that if $h_m\ne 0$, $m\in
\mathbb{N}_0$, then \textit{there always exists a unique set of the
numbers} $\lambda_k=\lambda_k(f,h,n)$ such that
$$
f(z)=\mathcal{H}_n(\{\lambda_k\},h;z) +O(z^n), \qquad z\to 0.
$$
On the other hand, in the above-mentioned papers $h$-sums are
also used as operators of differentiation, integration,
interpolation and extrapolation on certain classes of functions,
holomorphic in a fixed neighbourhood of the origin. In this case the
numbers $\lambda_k$ are already independent of individual functions
$f$ from the class and hence are of a universal kind. For instance, the following formulas for numerical differentiation
and integration, being exact for the polynomials of degree $\le
n-1$, are valid \cite{Dan2008}:
\begin{equation}
\label{300} zh'(z)\approx -h(z)+\sum_{k=1}^{n}\lambda_{1,k}
h(\lambda_{1,k} z);\quad \int_{0}^{z} h(t)\,dt\approx z
\sum_{k=1}^{n}\lambda_{2,k} h(\lambda_{2,k} z).
\end{equation}
Here the numbers $\lambda_{l,k}$ are absolute constants, being the
roots of the polynomials $P_{l,n}$ $(l=1,2)$, which can be defined
recursively as follows. Let $P_{l,0}=1$, $v_{l,1}=-1$ $(l=1,2)$,
then for $k=1,2,\ldots$ we have
$$
P_{l,k}=\lambda P_{l,k-1}+v_{l,k},\quad
v_{1,k}=-1-\sum_{j=1}^{k-1}\left(1-\frac{j}{k}\right)v_{1,j},\quad
v_{2,k}=-\frac{1}{k^2}-\sum_{j=1}^{k-1}\frac{v_{2,j}}{k(k-j)}.
$$

In 2013 we proposed \cite{CD-B,CD-K} a natural generalization of the $h$-sums, the so-called
 \textit{amplitude and frequency operators} (\textit{sums}) of the form
\begin{equation}
\label{gH} H_n(z)=H_n(\{\mu_k\},\{\lambda_k\},h;z):=\sum_{k=1}^n
\mu_k h(\lambda_k z),\qquad \mu_k,\lambda_k\in \mathbb{C},
\end{equation}
where \textit{amplitudes} $\mu_k$ and \textit{frequencies}
$\lambda_k$ are parameters, being independent of each other. In the
preprint \cite{CD-arxiv} and the present paper we give a detailed
exposition of the results announced in \cite{CD-B,CD-K}.

Later on, operators of the form (\ref{gH}) were studied in \cite{YF} but
with a fundamentally different approach for constructing them.
Namely, there was proposed not an analytic method for that as in
this paper, but a numerical one with small residuals (it will be
discussed in Section~\ref{Section6}).

The number $n$ in (\ref{gH}) is called \textit{the order} of the
amplitude and frequency operator, if there are no zeros among the
numbers $\mu_k$ and, moreover, the numbers $\lambda_k$ are pairwise
distinct (otherwise the order of the operator is less than $n$). As
in the case of $h$-sums, we regard amplitude and frequency sums  both as approximants of individual functions $f$, analytic
at the origin (and then $\mu_k=\mu_k(f,h,n)$ and
$\lambda_k=\lambda_k(f,h,n)$), and as special operators (of
differentiation, extrapolation, etc.), acting on certain classes of
functions (and then $\mu_k=\mu_k(n)$ and $\lambda_k=\lambda_k(n)$).

Introduction of the additional parameters $\mu_k$ enables us to
formulate the problem of $2n$-multiple (Pad\'{e}) interpolation at $z=0$
by means of the amplitude and frequency sums (in contrast to the $h$-sums, when
 $n$-multiple interpolation is only possible). Indeed, given
Maclaurin series
$$
f(z)=\sum_{m=0}^{\infty} f_m z^m, \qquad h(z)=\sum_{m=0}^{\infty}
h_m z^m, \qquad \text{where } f_{m} = 0, \text{ if } h_{m} = 0,
$$
we introduce the numbers $s_m=s_{m}(h,f)$:
\begin{equation}
\label{ssmm} s_{m}(h,f)=0, \hbox{ if } f_{m}=0; \qquad
s_{m}(h,f)=f_{m}/h_{m}, \hbox{ if } f_{m}\ne 0,\qquad m \in \mathbb{N}_0.
\end{equation}

For $|z|<\rho\cdot\min_{k=\overline{1,n}} |\lambda_k|^{-1}$ the
operator (\ref{gH}) has the form
\begin{equation*}
\label{gH+} H_n(z)=\sum_{k=1}^n \mu_k \sum_{m=0}^{\infty} h_m
(\lambda_k z)^m = \sum_{m=0}^{\infty}h_{m} \left(\sum_{k=1}^n
\mu_k \lambda_k^m\right) z^{m},
\end{equation*}
hence to realize the $2n$-multiple interpolation
\begin{equation}
f(z)=H_n(\{\mu_k\},\{\lambda_k\},h;z)+O(z^{2n}), \qquad z\to 0,
\label{2n-interpolation}
\end{equation}
or, which is the same,
\begin{equation}
f^{(m)}(z)=H_n(\{\mu_k\lambda_k^m\},\{\lambda_k\},h^{(m)};z)
+O(z^{2n-m}),\quad m=\overline{0, 2n-1}, \quad z\to 0,
 \label{2n-bis-interpol}
\end{equation}
the following conditions on the so-called \textit{generalized power
sums} (\textit{moments}) $S_m$ should be satisfied:
\begin{equation}
\label{SRS} S_m:=\sum_{k=1}^n \mu_k \lambda_k^m=s_m,\qquad
m=\overline{0,2n-1}.
\end{equation}

The system (\ref{SRS}) with unknown $\mu_k$, $\lambda_k$ and given
$s_m$ is well known as \textit{the discrete moment problem}.
Classical works of Prony, Sylvester, Ramanujan and papers of many
contemporary researchers  are devoted to the problem of its solvability
 (see \cite{Prony,Sylvester,Ramanujan,Lyubich,Lyubich2,Kung1}). Note that the
system (\ref{SRS}) is bound up with Hankel forms, orthogonal
polynomials, continued fractions, Gaussian quadratures and Pad\'{e}
approximants (a detailed review of these connections is given in
\cite{Lyubich,Lyubich2} and also in Section~\ref{par2}).

Suppose that the system (\ref{SRS}) is solvable. Then, following
\cite{Lyubich}, we call the system and its solution \textit{regular}
if all $\lambda_k$ are pairwise distinct and all $\mu_k$ are not
vanishing. In the case of regular systems (\ref{SRS}) we call the
problem of $2n$-multiple interpolation (\ref{2n-interpolation})
\textit{regularly solvable}. One of the methods for solving regular
systems~(\ref{SRS}) is due to Prony \cite{Prony}. Consider the
following product of determinants:
\begin{equation*}
\label{equ1} \left|
\begin{array}{ccccc}
1 & 0 & 0 & \ldots & 0\\
0 & \mu_1 & \mu_2 & \ldots & \mu_n\\
0 & \mu_1\lambda_1 & \mu_2\lambda_2 & \ldots & \mu_n\lambda_n\\
\ldots & \ldots & \ldots & \ldots & \ldots\\
0 & \mu_1\lambda_1^{n-1} & \mu_2\lambda_2^{n-1} & \ldots & \mu_n\lambda_n^{n-1}\\
\end{array}
\right| \cdot \left|
\begin{array}{ccccc}
1 & \lambda & \lambda^2 & \ldots & \lambda^n\\
1 & \lambda_1 & \lambda_1^2 & \ldots & \lambda_1^n\\
\ldots & \ldots & \ldots & \ldots & \ldots\\
1 & \lambda_n & \lambda_n^2 & \ldots & \lambda_n^n\\
\end{array}
\right|.
\end{equation*}
By regularity, the former of them does
not vanish and the latter does only for $\lambda=\lambda_k$ (as a
Vandermonde determinant). On the other hand, direct multiplication
of the determinants and
 taking into account (\ref{SRS}) give the following determinant, which is a polynomial of $\lambda$:
\begin{equation}
\label{G_n} G_n(\lambda):=\sum_{m=0}^n g_{m} \lambda^{m}= \left|
\begin{array}{ccccc}
1 & \lambda & \lambda^2 & \ldots & \lambda^n\\
s_0 & s_1 & s_2 & \ldots & s_n\\
s_1 & s_2 & s_3 & \ldots & s_{n+1}\\
\ldots & \ldots & \ldots & \ldots & \ldots\\
s_{n-1} & s_n & s_{n+1} & \ldots & s_{2n-1}\\
\end{array}
\right|.
\end{equation}
We call $G_n$ {\it a generating polynomial} (for functional
properties of such polynomials, including orthogonality,
completeness, etc., see \cite{Lyubich}). Consequently, the numbers
$\lambda_k$ are the simple roots of the generating polynomial $G_n$.
If the equation $G_n(\lambda)=0$ is solved, we substitute the
numbers found into the system (\ref{SRS}). Finally, extracting
any $n$ rows from (\ref{SRS}) leads to a linear system of equations
with unknowns $\mu_k$, which has a unique solution with
non-vanishing components. We now give known formulas for calculation
of the numbers $\mu_k$ (see, for example, \cite{DanDod2013}). Let
$\sigma_m$ and $\sigma_m^{(k)}$ denote elementary symmetric
polynomials of the form
$$
\sigma_m=\sigma_m(\lambda_1,\ldots,\lambda_n)= \sum_{1 \le
j_1<\ldots<j_m \le n}{\lambda_{j_1}\ldots \lambda_{j_m}}, \quad
m=\overline{1,n},
$$
$$
\sigma_0=1,\quad
\sigma_m^{(k)}=\sigma_m(\lambda_1,\ldots,\lambda_{k-1},0,
\lambda_{k+1},\ldots,\lambda_n), \quad k=\overline{1,n}.
$$
\begin{lemma}
\label{400} The numbers $\mu_k$ are the scalar products
$\mu_k=({\mathcal{L}}_k\cdot {\mathcal{S}})$, where ${\mathcal{S}}=
(s_0,\ldots, s_{n-1})$ and
\begin{equation}
\label{444} {\mathcal{L}}_k=\frac{{g}_n}{{G}'_n(\lambda_k)}
\left((-1)^{n-1}\sigma_{n-1}^{(k)},  \ldots ,
(-1)^{n-m}\sigma_{n-m}^{(k)}, \ldots , -\sigma_{1}^{(k)},\;
1\right).
\end{equation}
\end{lemma}
\begin{Proof}
If $V=V(\lambda_1,\lambda_2,\ldots,\lambda_n)$ is a Vandermonde
matrix of the first $n$ equations in the system (\ref{SRS}), then
 the elements of the $k$th row
${\mathcal{L}}_k$ of the matrix $V^{-1}$ have the form (\ref{444});
see \cite{DanDod2013} for more details.
\end{Proof}

We now formulate a known criterion of regularity in terms of
roots of the polynomial~$G_n$. Originally this criterion was obtained
in an algebraic form by Sylvester~\cite{Sylvester} (see also
\cite[Ch. 5]{Kung1}); later on Lyubich \cite{Lyubich,Lyubich2}
stated it in the analytical terms, which we use in the present paper.
\begin{thL}
\label{Criterion_Syl} The system $(\ref{SRS})$ is regular if and only if the generating polynomial $G_n$ is of degree $n$ and all
 its roots are pairwise distinct. Moreover, the regular system has a unique solution.
\end{thL}

This theorem immediately implies the following proposition about
 regular solvability of the interpolation problem (\ref{2n-interpolation}) for a function $f$, being analytic in a
neighbourhood of the origin.
\begin{theorem}
\label{th1} Suppose that the generating polynomial $G_n$,
constructed using the numbers $s_m=s_m(h,f)$, $m=\overline{0,2n-1}$,
is of degree $n$ and all its roots are pairwise distinct. Then the
amplitude and frequency operator $H_n$ is uniquely determined from the system
$(\ref{SRS})$  and realizes the $2n$-multiple interpolation
$(\ref{2n-interpolation})$ of the function~$f$ at the node~${z=0}$.
\end{theorem}
\begin{remark}
\label{remark_even_odd} If the function $f$ is even or odd, then the
local precision of the interpolation can be increased. Indeed, if $f$ is
even, then $f(z)=\tilde{f}(t)$, $t=z^2$, for some function
$\tilde{f}$, analytic at the point $t=0$, and the interpolation
 (\ref{2n-interpolation}) under the assumptions of Theorem~\ref{th1} with the basic function  $h$ gives
\begin{equation}
\label{interpolaton_even}
\tilde{f}(t)=\sum_{k=1}^n\mu_kh(\lambda_kt)+O(t^{2n})\quad
\Leftrightarrow \quad
f(z)=\sum_{k=1}^n\mu_kh(\lambda_kz^2)+O(z^{4n}),
\end{equation}
where it is necessary that $f_{2m}\neq 0 \Rightarrow h_m\neq 0$, see
(\ref{ssmm}). If $f$ is odd, then $f(z)=z\tilde{f}(t)$, $t=z^2$, and
analogously
$$
f(z)=\sum_{k=1}^n\mu_kzh(\lambda_kz^2)+O(z^{4n+1}).
$$

Note that, in contrast to similar discussions in \cite[Corollary~
2]{YF}, here we do not require the functions $f$ and $h$ to be of
the same parity.
\end{remark}
\begin{remark}
\label{remark2}
 One can consider the interpolation problem
(\ref{2n-interpolation}) with $O(z^{M})$, where $M>2n$ or $M<2n$, instead of $O(z^{2n})$. Then we accordingly
get overdetermined and underdetermined moment systems of the type
(\ref{SRS}):
\begin{equation}
\label{SRS_M} \sum_{k=1}^n \mu_k \lambda_k^m=s_m,\qquad
m=\overline{0,M-1}.
\end{equation}
In some cases the process of solving the consistent systems (\ref{SRS_M})
with $M\ne 2n$ can be reduced to the one of the standard
systems~(\ref{SRS}) (with $M=2n$). It can be done by eliminating the superfluous
equations or adding the missed ones (regarding this see
also \cite[\S5]{Lyubich}). But in the present paper we do not
consider the case $M\ne 2n$ for the following reasons.

The overdetermined systems (\ref{SRS_M}) belong to the non-regular
problems of the form (\ref{SRS}), where $2n$ is exchanged by $M$ and
$\mu_{k}=0$ for $k\ge n+1$. To apply the Prony-Sylvester method or
some other analytical approaches in the standard subsystem of the
system~(\ref{SRS_M}), one needs a preliminary analysis of its
consistency. As far as we know, there exist no reasonably general methods for this purpose. Moreover, it can be seen from the corresponding
overdetermined interpolation problem of the
form~(\ref{2n-interpolation}), where $O(z^{2n})$ is exchanged by
$O(z^{M})$, that its consistency is quite rigidly connected with
individual properties  of the functions $f$ and $h$ and thus one has
 a little chance to obtain more or less general interpolation formulas
of the overdetermined type.

For example, if a solvable system of the form (\ref{SRS}) is  supposed to
be consistent with the next equation $S_{2n}=s_{2n}$, then the
coefficient $f_{2n}$ cannot be chosen arbitrarily as it depends on the parameters $s_0,\ldots s_{2n-1},h_{2n}$ in a
certain way. Indeed, the following generalized Newton's formula,
connecting the coefficients $g_m$ of the generating polynomial $G_n$
and the moments $S_{v+m}$, is well-known:
\begin{equation}
\label{corollary_form_Newton++++} \sum_{m=0}^{n} S_{v+m}{g}_{m}=0,
\qquad v=0,1,\ldots.
\end{equation}
Therefore the values $s_m=s_m(h,f)$ (see (\ref{ssmm})) of the
moments $S_m$ with ${m>2n-1}$ (simultaneously with the coefficients
 $f_{m}=s_m h_m$) in a solvable overdetermined system of the form (\ref{SRS_M})
are uniquely determined from the system (\ref{SRS}) (we suppose that $h$ is fixed). A similar situation arises
 when one obtains formulas for numerical differentiation and
extrapolation (see the corresponding sections below). In these
problems the sequences $\{s_m\}$ have a certain arithmetic structure
and do not satisfy~(\ref{corollary_form_Newton++++}) for $v=n$ and $S_{2n}=s_{2n}$, i.e.
one gets solvable systems of the form (\ref{SRS}) with $M=2n$ although adding one more equation of the required form leads to an
inconsistent system (see, for example, (\ref{r_n_diff})).

As regards the underdetermined moment systems (when $M<2n$) and the
corresponding Pad\'{e} interpolation, they do not arise in the
present paper as we solve the interpolation problem of the order
higher than $M$ with the same number of independent parameters
 $\{\mu_k\}_{k=1}^n$ and $\{\lambda_k\}_{k=1}^n$.

Nevertheless, the systems (\ref{SRS_M}) (consistent and even
inconsistent) are of independent interest, moreover, they are actively studied in
numerical analysis. Various methods for finding an approximate  solution (in
different senses) to them  were developed, for example, in
\cite{Beylkin,Kung_residual,Beylkin2,Potts,YF} (characteristics
of one such method are discussed in
Section~\ref{Section6}).
\end{remark}

\section{Amplitude and frequency operators in classical problems}
\label{par2}

 We now consider several classical problems, which are
bound up with the class of amplitude and frequency operators.

\medskip

{\bf \thesection.1. Hamburger moment problem.} Theorem~\ref{th1}
raises the question about interpolating amplitude and frequency
operators with real $\lambda_k$ and $\mu_k$ (in particular, with
$\mu_k>0$). This question is well-studied \cite[Ch. 2]{Akhiezer1}
and can be settled by discretization of the following classical Hamburger
problem: given a sequence of real numbers $\{s_m\}$, $m\in
\mathbb{N}_0$, find a non-negative Borel measure $\mu$ on
$\mathbb{R}$ such that
\begin{equation}
\label{Hamburger} s_m = \int_{-\infty}^\infty \lambda^m\,d
\mu(\lambda),\qquad m\in \mathbb{N}_0.
\end{equation}
Namely, the following criterion is valid \cite[Ch. 2, \S1]{Akhiezer1}:
the problem (\ref{Hamburger}), where $m=\overline{0,2n-1}$, has a
unique solution with the spectrum, consisting of $n$ pairwise distinct
points $\lambda_1,\ldots,\lambda_n$, if and only if the leading
principal minors $\Delta_k$ of order $k$ of the infinite Hankel
matrix $(s_{i+j})_{i,j=0}^\infty$ satisfy the following conditions:
\begin{equation}
\label{Hamb_condition}
\Delta_1>0, \qquad \Delta_2>0,\qquad \ldots \qquad \Delta_{n}>0,\qquad \Delta_{n+1}=\Delta_{n+2}=\ldots=0.
\end{equation}
This implies that the discrete moment problem (\ref{SRS}) is regularly solvable in real numbers
$\lambda_k$ (this is equivalent to the fact that the polynomial
 (\ref{G_n}) has $n$ pairwise distinct real roots) and $\mu_k>0$ if and only if
 the sequence (\ref{ssmm}) satisfies the first $n$ inequalities in (\ref{Hamb_condition}).
Note that then the sequence $\{s_m\}_{m=0}^{2n-1}$ is called \textit{positive}.

\medskip

{\bf \thesection.2. Gauss and Chebyshev quadratures.} Given a
function $f$, analytic in a $\rho$-neighbourhood of the origin,
suppose that
$$
F(x):=\frac{1}{x}\int_{-x}^{x}f(t)\,dt,\qquad 0<x <\rho.
$$
To construct the amplitude and frequency operator
$H_n(\{\mu_k\},\{\lambda_k\},f;x)$ for $F(x)$, we
get the (positive) moment sequence $s_{m}=\frac{1-(-1)^{m+1}}{m+1}$,
$m=\overline{0,2n-1}$, from (\ref{ssmm})  and then consider the corresponding discrete
moment problem (\ref{SRS}). It is well-known that it is regular for
any $n$, moreover, the corresponding generating polynomial (\ref{G_n}) is the
Legendre polynomial $P_n$ (we write it in Rodrigues' form)
multiplied by a non-zero constant \cite[Ch. 7, \S2]{Krylov}:
 $$
G_n(x)=P_n(1) P_n(x), \qquad P_n(x):=\frac{1}{2^n
n!}\frac{d^n}{dx^n}(x^2-1)^n.
 $$
Therefore the frequencies $\lambda_k$ are real, pairwise distinct
and belonging to the interval $(-1,1)$ (as roots of the Legendre
polynomials, forming an orthogonal system on $[-1,1]$). The
amplitudes $\mu_k$ are determined via the numbers $\lambda_k$ by
the well-known formulas \cite[Ch. 10, \S 3]{Krylov}:
$$
\mu_k = \frac{2}{\left( 1-\lambda_k^2 \right)
[P'_n(\lambda_k)]^2}>0.
$$
Thus we obtain the interpolation formula
\begin{equation}
\label{GAUSS}
\frac{1}{x}\int_{-x}^{x}f(t)\,dt=\sum_{k=1}^n\mu_kf(\lambda_kx)+r_n(x),
\qquad r_n(x)=O(x^{2n}),
\end{equation}
which is a Gaussian quadrature for each fixed $x$. The amplitudes
and frequencies depend only on $n$ but not on $f$. It is known
\cite[Ch. 10]{Krylov} that the Gaussian quadratures are of the
highest algebraic degree of accuracy among all the formulas of the
form (\ref{GAUSS}) and exact (i.e. $r_n(x)\equiv 0$) for the
polynomials of degree $\le 2n-1$.

In a similar manner one can obtain interpolation formulas for
integrals with classical weights. For example, for
$$
F(x):=\int_{-x}^{x}\frac{f(t)}{\sqrt{x^2-t^2}}\,dt,\qquad 0< x
<\rho,
$$
we have
$$
s_{2m}=\pi\frac{(2m)!}{(2^m m!)^2},\quad s_{2m+1}=0,\quad
m=\overline{0,n-1},\qquad G_n(x)=(-2^{1-n}\pi)^n T_n(x),
$$
where $T_n(x)=\cos(n\arccos x)$ are the Chebyshev polynomials of the
first kind. Calculating the amplitudes $\mu_k$ via the
frequencies $\lambda_k$ leads to the following Gauss-Chebyshev quadrature \cite[\S 25.4.38]{Abramovitz} for real
$x$, $0<x <\rho$:
\begin{equation}
\label{Gauss-Chebyshev}
\int_{-x}^{x}\frac{f(t)}{\sqrt{x^2-t^2}}\,dt=\frac{\pi}{n}
\sum_{k=1}^n f(\lambda_kx)+r_n(x), \quad \lambda_k=\cos
\tfrac{2k-1}{2n}\pi, \quad r_n(x)=O(x^{2n}),
\end{equation}
whose characteristic property is the equality of the amplitudes
$\mu_k=\pi/n$. The remainder can be written more precisely:
\begin{equation}
\label{Gauss-Chebyshev-error} r_n(x)=\frac{\pi
f^{(2n)}(\xi)}{2^{2n-1}(2n)!}x^{2n},\qquad \xi\in(-x,x).
\end{equation}
One can deduce it from \cite[\S 25.4.38]{Abramovitz} by a suitable
change of variables.

\medskip

{\bf \thesection.3. Pad\'{e} approximants.} The Pad\'{e}
approximants as well as the Gaussian quadratures are closely related
to the classical moment problem (see, for instance, \cite{Dzyadyk}).
Construction of the amplitude and frequency sum
$H_n(\{\mu_k\},\{\lambda_k\},h;x)$ with the basis function
$h(z)=(z-1)^{-1}$ for some function $f$, analytic in a neighbourhood of the
origin, leads to the sequence of moments $s_m=-f_m$,
$m=\overline{0,2n-1}$. If the generating polynomial (\ref{G_n}) for
this sequence is of degree $n$ and all its roots are pairwise
distinct, then by Theorem~\ref{th1} we get the following
interpolation identity:
\begin{equation}
\label{2}
f(z)=\sum_{k=1}^n\frac{\mu_k}{\lambda_kz-1}+O(z^{2n}).
\end{equation}
This is a classical Pad\'{e} approximant of order $[(n-1)/n]$. We
recall that classical Pad\'{e} approximants of order $[m/n]$ are
interpolating rational functions of the form $P_m(f;z)/Q_n(f;z)$
(see, for instance, \cite[\S1.1]{Baker}). Note that the method for
solving the problem (\ref{SRS}), proposed by Ramanujan
\cite{Ramanujan}, is equivalent to the one for construction of the
interpolation formula (\ref{2}) (see \cite{Lyubich, Lyubich2}).

\medskip

{\bf \thesection.4. Exponential sums.} Let
$h(z)=\exp(z)$ be a  basis function  in the amplitude and frequency operator $H_n(\{\mu_k\},\{\lambda_k\},h;z)$ and $f$ be a function, which we are going to interpolate.
 The corresponding sequence of moments is $s_m=m!f_m$,
$m=\overline{0,2n-1}$. Suppose that the problem (\ref{SRS}) for this
sequence is regular. Then the following formula of $2n$-multiple
interpolation at the origin  holds:
$$
f(z)=\sum_{k=1}^n \mu_k e^{\lambda_kz}+O(z^{2n}).
$$
(In particular, this result has been already obtained
in \cite{Buchmann} and \cite{Lyubich}.)

Interpolation of functions by exponential sums with \textit{simple}
equidistant nodes was considered by Prony~\cite{Prony}. At present,
many works are devoted to this method and its various modifications
and applications  (see, for instance,
\cite{Korobov,Beylkin2,Beylkin}, \cite[Ch. 6]{Braess} and references there).

A vast investigation of  \textit{the exponential series} was
conducted in the scientific school of Leont'ev \cite{Leontiev2}. It
is worth mentioning here that members of the school also actively
studied several \textit{generalizations of the exponential series}
 (see, for instance, \cite{Gromov2,Shevtsov,Leontiev3}). Namely, they
enquired into the problem of completeness of the infinite systems
$\{h(\lambda_k z)\}$, where $h$ were entire functions and
$\lambda_k$ given numbers. Consequently, they actually considered
some \textit{representations} of analytic functions  $f$ by amplitude
and frequency sums of infinite order,
$H_\infty(\{\mu_k\},\{\lambda_k\},h;z)$, and properties of these
representations (domains of convergence, admissible classes of the
numbers $\lambda_k$ and functions $f$, connections between $\mu_k$
and $\lambda_k$, etc.). In contrast to this approach, we consider
\textit{approximations} by amplitude and frequency sums of
finite order and respective errors. Moreover, the parameters $\mu_k$
and $\lambda_k$ are not given but uniquely determined by the
functions $f$ and $h$. Furthermore, in different applications we
regard amplitude and frequency sums
 as operators with fixed (universal) numbers  $\mu_k$ and $\lambda_k$, being determined by the analytic nature of these operators.

 \section{Examples}
In this section we give several examples of approximating amplitude
and frequency sums for some special functions, in particular, Bessel
functions (all arising discrete moment problems are regularly
solvable). We will also compare our approximants
with similar ones, obtained by other authors.
\begin{example}
It is known \cite[\S 9.1.20-21]{Abramovitz} that the Bessel
function of order zero, $J_0$, has the representation
$$
J_0(\pm x)=\frac{1}{\pi}\int_{-x}^{x}
\frac{\exp(it)}{\sqrt{x^2-t^2}}\,dt=
\frac{1}{\pi}\int_{-x}^{x}\frac{\cos(t)}{\sqrt{x^2-t^2}}\,dt, \qquad
x>0.
$$
From this, (\ref{Gauss-Chebyshev}) and (\ref{Gauss-Chebyshev-error})
we get the following  high-accuracy local approximation of $J_0$ at
the point $x=0$ by an amplitude and frequency sum of order~$n$:
\begin{equation}
\label{J_0} J_0(x)=\frac{1}{n}\sum_{k=1}^n \cos(x\cdot\cos
\tfrac{(2k-1)\pi }{2n})+r_n(x), \qquad |r_n(x)|\le
\frac{|x|^{2n}}{2^{2n-1}(2n)!}.
\end{equation}
Note that one can obtain the same formula by interpolation of the series
\cite[\S 9.1.10]{Abramovitz}
\begin{equation}
\label{J_0000}
J_0(z)=\sum_{m=0}^\infty\frac{(-1)^m}{(2^mm!)^2}z^{2m}
\end{equation}
by amplitude and frequency sums with the basis function
 $h(z)=\exp(z)$. Furthermore, then
 $$
 s_{2m}=\frac{(2m)!}{(2^m m!)^2},\qquad s_{2m+1}=0,\qquad
 m=\overline{1,n-1},
 $$
cf. $\{s_{2m}\}$  for the Gauss-Chebyshev quadrature in Section~2.2.

If in (\ref{J_0}) we use $2n$ instead of $n$ and take into account
the symmetry of the frequencies obtained and the parity of cosine,
then we get the following amplitude and frequency sum (of the same
order $n$):
\begin{equation}
\label{J_0_super} H_n(x)=\frac{1}{n}\sum_{k=1}^{n} \cos(x\cdot\cos
\tfrac{(2k-1)\pi}{4n}),
\end{equation}
for which by (\ref{Gauss-Chebyshev-error}) we have
\begin{equation}
\label{J_0_super_error} J_0(x)\approx H_n(x),\qquad
|J_0(x)-H_n(x)|\le \frac{|x|^{4n}}{2^{4n-1}(4n)!}.
\end{equation}
Note that the formula (\ref{J_0_super}) can be also obtained via
(\ref{interpolaton_even}).

In \cite[Table 1 and Formula (25)]{YF} the following approximant was
obtained (for the approach from \cite{YF} see also
Section~\ref{Section6}):
\begin{equation}
\label{YF_J_0_cos} J_0(x)\approx
\omega_{11}(x):=\sum_{k=1}^{11}\alpha_m \cos(\beta_m x), \quad
\alpha_k\in(0.086,0.096),\quad \beta_k\in(0,1).
\end{equation}
In comparison with this approximant, the sum (\ref{J_0_super})  gives much more
precise results for $n=11$ in the segment $[0,30]$. Calculations in
 Maple show that
\begin{equation*}
\label{YF_J_0_cos++}
 {\rm log}_{10}\frac{|J_0(x)-\omega_{11}(x)|}{|J_0(x)-H_{11}(x)|}>
 M(x):=\frac{38}{x+1},
\qquad x\in [0;30].
\end{equation*}
Moreover, in the segment [0,10] the minorant $M(x)$ can be exchanged
by the more precise $\tilde M(x)=44\log_{10}
\frac{1}{x}+50\to\infty$, $x\to 0$. Note that the absolute error of
the formula (\ref{YF_J_0_cos}) in the segment $[0,10]$ is quite
small and close to $\varepsilon=10^{-12}$.

For $30\le x\le 40$ the absolute error of the formula
 (\ref{J_0_super}) is also less than the one of (\ref{YF_J_0_cos}), but for $x\ge 40$ the errors of both approximants can reach  $10^{-1}$,
 thus use of them does not seem to be fair for such $x$.
\end{example}
\begin{example}
From (\ref{2n-bis-interpol}) and (\ref{J_0_super}) we get the
following interpolation formula for the derivative of the Bessel
function, $J'_0$ (see \cite[\S 9.1.28]{Abramovitz}):
\begin{equation}
\label{derivative_J_0++} J'_0(x)\approx H'_n(x)=
-\frac{1}{n}\sum_{k=1}^{n}
 \cos \tfrac{(2k-1)\pi}{4n}\, \sin(x\cdot\cos \tfrac{(2k-1)\pi}{4n})
\end{equation}
with the error (see (\ref{Gauss-Chebyshev-error}))
\begin{equation}
\label{derivative_J_0_super_error} |J'_0(x)-H'_n(x)|\le
\frac{|x|^{4n-1}}{2^{4n-1}(4n-1)!}.
\end{equation}

For comparison, look at the approximant from \cite[Table 2 and Formula (30)]{YF}:
\begin{equation}
\label{YF_J_1}
 J'_0(x)\approx \Omega_{13}(x):=-\sum_{k=1}^{13}\alpha_m \sin(\beta_m
x),\quad \alpha_k\in(0.001,0.08),\quad \beta_k\in(0.12,1).
\end{equation}
Calculations show that the formula (\ref{derivative_J_0++}) with the
same $n=13$ has noticeably less error in the segment $[0,40]$:
\begin{equation*}
\label{YF_J_0_cos++++}
 {\rm log}_{10}\frac{|J'_0(x)-\Omega_{13}(x)|}{|J'_0(x)-H'_{13}(x)|}>
 M_1(x):=\frac{40}{x+1},
\qquad x\in [0,40].
\end{equation*}
Moreover, in the segment $[0,10]$ the minorant $M_1(x)$ can be
exchanged by the more precise $ \tilde M_1(x)=43\log_{10}
\frac{1}{x}+48\to\infty$, $x\to 0$. Note that the absolute error of
(\ref{YF_J_1}) in $[0,10]$ is close to $\varepsilon=10^{-12}$.

For $40\le x\le 45$ the absolute error of the formula
 (\ref{derivative_J_0++}) is also less than the one of (\ref{YF_J_1}). For $x\ge 45$ the absolute error of both formulas
 can be $10^{-2}$ and thus use of them may be unreasonable.
\end{example}
\begin{example}
Let us obtain one more representation for the function $J_0$ by taking
\begin{equation} \label{sinc}
\textrm{sinc}\; x=\frac{\sin x}{x}=\sum_{m=0}^\infty
\frac{(-1)^m}{(2m+1)!}x^{2m}
\end{equation}
as a basis function and using the approach for interpolation of even
functions from Remark~\ref{remark_even_odd}. Namely, we take into
account that $J_0$ and $\textrm{sinc}$ are even and interpolate
the function $\tilde{f}(t):=J_0(x)$, $t=x^2$, by amplitude and
frequency sums with the basis function ${h}(t)=\textrm{sinc}\; x$.
As it can be easily checked (see (\ref{J_0000}) and (\ref{sinc})),
then $s_m=\frac{(2m+1)!}{(2^mm!)^2}$, $m=\overline{0,2n-1}$. Solving
the corresponding problem (\ref{SRS}) (in all examples considered we obtained non-negative $\lambda_k$ and real
$\mu_k$) yields
$$
\tilde{f}(t)=\sum_{k=1}^n\mu_k h(\lambda_k t)+r_n(t),\qquad
r_n(t)=O(t^{2n}),
$$
or, which is the same,
\begin{equation}
\label{J_0_sinc} J_0(x)=\sum_{k=1}^n\mu_k \,
\textrm{sinc}(\sqrt{\lambda_k}x)+r_n(x),\qquad r_n(x)=O(x^{4n}).
\end{equation}

We now compare (\ref{J_0_sinc}) with the approximate equality
 $J_0(x)\approx
\sum_{k=1}^{11}\alpha_k\,\textrm{sinc}(\beta_k z)$, obtained in
\cite[Table 3 and Formula (32)]{YF} (the amplitudes and frequencies
were found there by a numerical method). It turns out that the
amplitude and frequency sum in
 (\ref{J_0_sinc}) with the
same order $n=11$  gives more precise results for $x\in[0;42]$
(especially in a small neighbourhood of $x=0$). For $x\ge 40$ the
absolute errors of both formulas can already exceed  $10^{-3}$.
\end{example}
\begin{remark}
The authors of \cite{Beylkin,Potts} consider numerical methods for
interpolation of the function $J_0$ in equidistant nodes by
amplitude and frequency sums of the form $\sum_{k=1}^{n}\alpha_k
\exp(\beta_k x)$ with complex (but not purely imaginary) $\alpha_k$
and $\beta_k$ and quite fast-decreasing moduli of exponents. This
enables them to get good approximants of $J_0$ on large intervals.
For example, in \cite{Beylkin} such an approximant  with $n=28$ has
the error $\varepsilon\le 10^{-10}$ for $x\in[0,100\pi]$. The
approximant with $n=20$,  obtained in \cite[Example~4.5]{Potts} by
slightly different methods, has the error $\varepsilon\le 10^{-4}$
for $x\in [0,1000]$.

Note that the sums (\ref{J_0_super}) and (\ref{derivative_J_0++})
are not recommended for use on such large intervals because of their
local character. To guarantee a reasonable rate of approximation,
they must have the order comparable with length of the intervals considered
 (it can be seen from the order-precise estimates  (\ref{J_0_super_error})
 and (\ref{derivative_J_0_super_error})). For example, for $n=20$ they do not give
 reasonable quality of approximation in the segment $[0,1000]$ but in the subsegments $[0,40]$, $[0,30]$ and $[0,20]$ the approximation errors
 are less than $10^{-16}$, $10^{-25}$ and $10^{-38}$, correspondingly.
\end{remark}

\section{Analytic regularization of the interpolation by amplitude and frequency operators}
\label{regularization_section}
 {\bf \thesection.1. Variation of the moments.} The $2n$-multiple interpolation problem has great
difficulty in the case when the conditions of regularity (see
Theorem~\ref{th1}) are not satisfied. In particular, it can
be inconsistent then. In order to avoid this difficulty, we propose
a method for analytic regularization of the discrete moment problem.
It consists in a certain variation of the right hand sides $s_m$ of
the system (\ref{SRS}), namely, in adding the generalized
power sums ${\sigma}\sum_{k=1}^{\nu}\alpha_k (r\beta_k)^m$ to them
(another approach is described in Remark~\ref{REG+++}). The
parameters $\alpha_k$, $\beta_k$ are independent of  $s_m$, and
$\sigma$ and $r$ depend only on $\max \{|s_m|\}$.  From the point of
view of the interpolation problem, this corresponds to the fact that
we get a new regularly solvable problem of the form
(\ref{2n-interpolation}) (and (\ref{SRS})), where $f$ is exchanged
for a new varied function $\tilde f$ such that  $\tilde f(z)-f(z)$ is
the amplitude and frequency sum $\sigma H_{\nu}(\{\alpha_k\},\{r
\beta_k\},h;z)$. We emphasize that  the above-mentioned variation of
$s_m=s_{m}(h,f)$ in the framework of the interpolation problem is
universal as it depends only on $\max \{|s_m|\}$ but
not on the functional properties of  $f$ and $h$. Moreover, as we
will see below, if $\alpha_k$, $\beta_k$, $\sigma$, $r$ are chosen
appropriately, then the difference $\tilde f-f$ is just a certain
binomial. Let us describe the regularization method in detail.
Instead of the function $f$, whose corresponding problem (\ref{SRS}) is
not regular, we introduce the varied function
\begin{equation}
\label{F(z)} {\tilde f}(z):=f(z)+\sigma H_{\nu}(\{\alpha_k\},\{r
\beta_k\},h;z),\qquad p\in \mathbb{C},\qquad {\nu}\in \mathbb{N},
\end{equation}
where $\alpha_k$, $\beta_k$ $\sigma$, $r$ are constants. Since
$$
{\tilde f}(z)= \sum_{m=0}^\infty f_m z^m
+\sigma\sum_{k=1}^{\nu}\alpha_k \sum_{m=0}^\infty h_m (r
\beta_kz)^m= \sum_{m=0}^\infty \left(s_m+\tau_m \right)h_m z^m,
$$
where $\tau_m:=\sigma\sum_{k=1}^{\nu}\alpha_k (r\beta_k)^m$ and
 $s_m=s_{m}(h,f)$ as above (see (\ref{ssmm})), instead of  $(\ref{SRS})$ we obtain the system
\begin{equation}
\label{Var_disc_mom_syst} \sum_{k=1}^n \mu_k
\lambda_k^m=s_m+\tau_m, \qquad m=\overline{0,2n-1},
\end{equation}
which differs from the system  (\ref{SRS}) in the regularizing
summands $\tau_m$ (if the initial system is regular, then it is
natural to set $\sigma=\tau_m\equiv 0$). The task is now to find
the numbers $\tau_m$ such that the conditions of Theorem~\ref{th1}
are satisfied. Assume that we have done this, then by
Theorem~\ref{th1} we obtain the interpolation identity
$$
{\tilde f}(z)=H_{n}(\{\mu_k\},\{\lambda_k\},h;z)+O(z^{2n}).
$$
Returning to the initial interpolation problem and taking into account (\ref{F(z)}) yield
\begin{equation}
\label{f_reg} f(z)=H_n(\{\mu_k\},\{\lambda_k\},h;z)-\sigma
H_{\nu}(\{\alpha_k\},\{r\beta_k\},h;z) +O(z^{2n}).
\end{equation}
At the same time it is reasonable to choose ${\nu}$ as small as
possible. However, there is a natural restriction on $\nu$, which is
seen from the following statement (cf. \S4 in \cite{Lyubich}).
\begin{lemma}
\label{ranks} It is necessary for the regularity of the varied system $(\ref{Var_disc_mom_syst})$ that
$$
{\nu}\ge n-\rank \left(s_{i+j-2}\right)_{i,j=1}^n.
$$
\end{lemma}
\begin{Proof}
Consider the Hankel matrices
$$
\mathbf{H}:=\left(s_{i+j-2}\right)_{i,j=1}^n, \qquad
\mathbf{R}:=\sigma\left(\sum_{k=1}^{\nu} \alpha_k (r
\beta_k)^{i+j-2}\right)_{i,j=1}^n.
$$
For the regularity of the system $(\ref{Var_disc_mom_syst})$ it is necessary that the coefficient before $\lambda^n$
 of the corresponding generating polynomial $G_n$ does not vanish,
i.e. $\det (\mathbf{H}+\mathbf{R})\ne 0$, $\rank
(\mathbf{H}+\mathbf{R})= n$. It is well known \cite[\S 0.4.5]{Horn}
that for any $n\times n$-matrices $A$ and $B$
\begin{equation}
\label{rank_inequality} \rank (A+B)\le \min\{\rank A + \rank B;n\}.
\end{equation}
Consequently, if the system $(\ref{Var_disc_mom_syst})$ is regular,
then necessarily $n \le \rank \mathbf{H} + \rank \mathbf{R}$. It
remains to note that $\rank \mathbf{R}\le {\nu}$. Indeed, the
following representation is valid:
$$
\mathbf{R}=\sigma\sum_{k=1}^{\nu} \alpha_k \mathbf{C}(k), \qquad
\mathbf{C}(k):=\left((r\beta_k)^{i+j-2}\right)_{i,j=1}^n,
$$
where $\rank \mathbf{C}(k)=1$ (each next row is the previous one
multiplied by $r\beta_k$). From this by the property
(\ref{rank_inequality}) we obtain the required bound for $\rank
\mathbf{R}$.
\end{Proof}

\medskip

{\bf \thesection.2. Parameters of the regularization.} In the problems under consideration the
ranks of matrices $\left(s_{i+j-2}\right)_{i,j=1}^n$ are small.
Consequently, taking into account Lemma~\ref{ranks},  we will
consider only the overall case ${\nu}=n$ to solve the regularization
problem. Then the formula (\ref{f_reg}) has $2n$ summands (if
$\sigma\ne 0$, $r\ne 0$), and in this sense the amplitude and
frequency sums obtained  have no advantage over the $h$-sums of
order $2n$. However, an appropriate choice of $\alpha_k$, $\beta_k$,
$\sigma$ and $r$ can essentially simplify the latter sum in
(\ref{f_reg}). Indeed, let $p$ and $q$ be fixed non-zero complex
numbers and
\begin{equation}
\label{alpha-beta} \alpha_k=\beta_k=\exp \left(\frac{2\pi (k-1)
i}{n}\right),\qquad r=\left(\frac{q}{p}\right)^{1/n},\qquad
\sigma=\frac{p^2}{nq} r,\qquad k=\overline{1,n},
\end{equation}
where the number $r$ is any of the $n$ values of the root. Then, as it can
be easily seen,  in (\ref{Var_disc_mom_syst}) we obtain
\begin{equation}
\label{TAU} \tau_{n-1}=p,\quad \tau_{2n-1}=q;\qquad \tau_{m}=0 \quad
\hbox{   for other   }\quad  m.
\end{equation}
Indeed,
$$
\tau_m=\frac{p^2}{nq} r^{m+1}\sum_{k=1}^n \exp\left({\frac{2\pi
(k-1) i}{n}(m+1)}\right),
$$
where the sum of the exponents equals $n$ if  $m+1$  is
divisible by $n$ and zero if not. Consequently,
\begin{equation*}
\label{SVORACH} \sigma H_n(\{\alpha_k\},\{r\beta_k\},h;z)= p\,
h_{n-1}z^{n-1}+q\,h_{2n-1}z^{2n-1}+O(z^{2n}).
\end{equation*}
Thus, assuming the regularity of the varied problem (\ref{Var_disc_mom_syst}), we get the formula
\begin{equation}
\label{f_reg_expansion} f(z)=H_n(\{\mu_k\},\{\lambda_k\},h;z)- p\,
h_{n-1}z^{n-1}-q\,h_{2n-1}z^{2n-1}+O(z^{2n}).
\end{equation}

In order to obtain the main result of this section, we now show that the
above-mentioned problem is indeed regular for a certain choice of
the parameters $p$  and $q$. Below we give a possible way of such a
choice.  The generating polynomial of the system
(\ref{Var_disc_mom_syst}) for $\alpha_k$ and $\beta_k$ from
(\ref{alpha-beta}) has the form
\begin{equation}
\label{modif_equ} {G}_n(\lambda)={G}_n(p,q;\lambda):= \left|
  \begin{array}{ccccc}
  1 & \lambda & \cdots & \lambda^{n-1} & \lambda^n \\
    s_0 & s_1 & \cdots & s_{n-1}+p & s_{n} \\
    s_1 & s_2 & \cdots & s_n & s_{n+1} \\
 \cdots & \cdots & \cdots & \cdots & \cdots \\
   s_{n-1}+p & s_{n} & \cdots & s_{2n-2} & s_{2n-1}+q \\
  \end{array}
\right|.
\end{equation}
Obviously, for the parameters $p$ and $q$, being sufficiently large
in modulus (comparing with the moments
 $s_k$ and independently of each other), the roots of this polynomial are arbitrarily close to
 those of the polynomial of the form (\ref{modif_equ}) with $s_k=0$, $k=\overline{0,2n-1}$. The latter polynomial,
  as it can be easily checked by expanding the determinant along the first row, has the form
$$
(-1)^{n(n+1)/2}p^{n}(\lambda^n-q/p),
$$
and all its $n$ roots are pairwise distinct. From here, the formula
(\ref{f_reg_expansion}) and Theorem~\ref{th1} we obtain the following result.
\begin{theorem}
\label{th2} Given $p$ and $q$ sufficiently large in modulus, the
varied problem $(\ref{Var_disc_mom_syst})$ with the parameters
$(\ref{TAU})$ has a regular solution $\{\mu_k\}$, $\{\lambda_k\}$.
Moreover, for the constants $c_1=-ph_{n-1}$ and $c_2=-qh_{2n-1}$ the
following interpolation formula holds:
$$
f(z)=c_1z^{n-1}+c_2 z^{2n-1}+\sum_{k=1}^n\mu_k h(\lambda_k
z)+O(z^{2n}).
$$
\end{theorem}
\begin{remark}
\label{REG+++} The above-mentioned regularization with the
parameters from (\ref{alpha-beta}) and (\ref{TAU}) is actually
equivalent to adding the binomial $c_1z^{n-1}+c_2 z^{2n-1}$
with non-vanishing coefficients $c_{1}$ and $c_{2}$ to the function
$f$. In what follows we will expand the class of
regularizable problems by showing that $c_1$ and $c_2$ can be chosen in a different way and not necessarily non-vanishing. In particular, in
the extrapolation problem from Section~\ref{par-extrap} it will be reasonable to set $c_2=0$.
\end{remark}
\begin{remark}
\label{rem2} The conditions on $p$ and $q$, mentioned in
Theorem~\ref{th2}, are quite qualitative and need additional
specification in practice.
 Several methods for this will be proposed below in particular applications. In the general case one can use the following observations. The leading coefficient $g_n=g_n(p)$ of the polynomial ${G}_n$ is obviously a polynomial of $p$ of degree $n$, hence $\deg
{G}_n(p,q;\lambda)=n$ for all $p$ except those from the set
\begin{equation}
\label{Pi-general} \Pi:=\{p: g_n(p)=0\},
 \end{equation}
containing no more than  $n$ points. It is possible to obtain some estimates for the boundaries of
the set $\Pi$ using that a matrix with strict diagonal dominance is non-singular
 (see the Levy–Desplanque theorem in \cite[Th. 6.1.10]{Horn}). Namely, if we choose
$p$ satisfying the inequality
$$
|s_{n-1}+p|>\sum_{j=1, j\neq n-i+1}^{n}|s_{i+j-2}|, \qquad
i=\overline{1,n},
$$
then the determinant for $g_n$ is strict
diagonally dominant and hence $g_n\ne 0$. For this it is sufficient to take,
for example,
$$
 p> n \max_{k=\overline{0,2n-1}} |s_k|.
$$
We now suppose that the generating polynomial (\ref{modif_equ}) is of degree $n$. Then the question about ``separation'' of its multiple roots arises. As it is easily seen,
$$
{G}_n(p,q;\lambda)={\mathcal S}(p;\lambda)
+q\mathcal{T}(p;\lambda),
$$
where the polynomial ${\mathcal S}$ is of degree $n$ and the
polynomial ${\mathcal T}$ is of degree $\le n-1$; both polynomials
depend only on  $p$. The following statement is valid.
\begin{lemma}
\label{RAZDEL} Suppose that in each multiple root $($if any$)$ of the polynomial ${G}_n$  the polynomial ${\mathcal T}$ either does not vanish or has a simple root. Then there exists an arbitrarily small variation
$\delta\ne 0$ of the parameter  $q$ such that the polynomial ${G}_n(p,q+\delta;\lambda)$ has  $n$ simple roots.
\end{lemma}
\begin{Proof}
Assume that ${\lambda}_0$ is an $s$-multiple ($s\ge 2$) root of the
polynomial ${G}_n(p,q;\lambda)$. Then in a sufficiently small
neighbourhood of the root the polynomial
$$
{G}_n(p,q+\delta;\lambda)={G}_n(p,q;\lambda) +\delta{\mathcal
T}(p;\lambda)
$$
 has the form
$$
{G}_n(p,q+\delta;\lambda)=({\lambda}-{\lambda}_0)^s({\alpha}+
O({\lambda}-{\lambda}_0))+\delta
(t_0+t_1(\lambda-{\lambda}_0)+O(({\lambda}-{\lambda}_0)^2)),
$$
where ${\alpha}\ne 0$, $|t_0|+|t_1|\ne 0$ and values
$O({\lambda}-{\lambda}_0)$, $O(({\lambda}-{\lambda}_0)^2)$,
$\lambda\to \lambda_0$, are independent of $\delta$. Choose small
$\varepsilon>0$ and $\delta=\delta(\varepsilon)$ so that in the disc
${|\lambda-\lambda_0|\le 2\varepsilon}$ the polynomial ${G}_n$ has
no roots, distinct from  $\lambda_0$ (we take into account that the
roots  depend on $\delta$ continuously), and
$$
|{G}_n(p,q;\lambda)|>|\delta{\mathcal T}(p;\lambda)|,\qquad
|\lambda-\lambda_0|=\varepsilon.
$$
By Rouch\'{e}’s theorem, the polynomial ${G}_n(p,q+\delta;\lambda)$
has exactly $s$ roots in the disc
${|{\lambda}-{\lambda}_0|<\varepsilon}$; we will use $\tilde
\lambda_k$ to denote them. If $t_0\ne 0$, then these roots satisfy
the equation
$$
({\lambda}-{\lambda}_0)^s=-\frac{\delta t_0}{\alpha}\;
(1+O(\varepsilon)),\qquad \varepsilon\to 0.
$$
If $t_0=0$, $t_1\ne 0$, then $\tilde \lambda_1=\lambda_0$ and other roots satisfy the equation
$$
({\lambda}-{\lambda}_0)^{s-1}=-\frac{\delta t_1}{\alpha}\;
(1+O(\varepsilon)),\qquad \varepsilon\to 0.
$$
In any case we get $s$ simple roots. Suppose that $\varepsilon$ and
$|\delta|$ are so small that the above-mentioned method works
simultaneously for all the multiple roots but all the simple ones remain
simple (it is possible as the roots  depend on $\delta$
continuously). Then we get the polynomial
${G}_n(p,q+\delta;\lambda)$ with $n$ simple roots.
\end{Proof}
\end{remark}
It follows from the aforesaid  that the following conjecture is very
likely: the set of the parameters $(p,q)$, for which the
interpolation problem considered in this section is regularly
solvable, is everywhere dense in $\mathbb{C}^2$. But now we have only the following statement, which is a supplement to
Theorem~\ref{th2}.
\begin{theorem} \label{th22}
Suppose that $p\notin \Pi$ $($see $(\ref{Pi-general})$$)$ and the
conditions of Lemma~$\ref{RAZDEL}$ are satisfied. Then there exists
an arbitrarily small variation  $\delta\ne 0$ of the parameter $q$
such that the varied problem  $(\ref{Var_disc_mom_syst})$ with
$\tau_{n-1}=p$, $\tau_{2n-1}=q+\delta$ $($all other $\tau_{m}=0)$
has a regular solution $\{\mu_k\}$, $\{\lambda_k\}$, and for the
constants $c_1=-ph_{n-1}$ and $c_2=-(q+\delta)h_{2n-1}$ the
following interpolation formula holds:
$$
f(z)=c_1z^{n-1}+c_2 z^{2n-1}+\sum_{k=1}^n\mu_k h(\lambda_k
z)+O(z^{2n}).
$$
\end{theorem}

\section{Numerical differentiation by amplitude and frequency operators}
\label{par-diff}
 {\bf \thesection.1. Statement of the problem.} As an application of the
regularization method we consider the problem of $2n$-multiple
interpolation of the function  $zf'(z)$ by amplitude and frequency
operators $H_n$ with the basis function $f$.
 (As above we suppose that $f$ is defined and holomorphic in a neighbourhood
 of the origin.) The solution to this problem would allow us to obtain a high-accuracy formula
 for
 numerical differentiation with local precision $O(z^{2n})$.
 However, the discrete moment problem (\ref{SRS}) with $s_m=m$, $m=\overline{0,2n-1}$, which we get in this case, is non-regular
(the generating polynomial (\ref{G_n}) is of the degree less than $n$
for $n=1$ and $n\ge 3$ as the algebraic adjunct
 to $\lambda^n$ obviously vanishes, and has the double root $\lambda=1$ for $n=2$; both cases do not satisfy Theorem~\ref{th1}).
Here we apply the regularization method mentioned in
Remark~\ref{REG+++}. More precisely, given some complex parameters
$p$ and $q$, we consider the varied function
\begin{equation}
\label{TILDE-f} \tilde{f}(z):=zf'(z)+ p\,
f_{n-1}z^{n-1}+q\,f_{2n-1}z^{2n-1}, \qquad
zf'(z)=\sum_{m=0}^{\infty} mf_m z^m
\end{equation}
and the interpolating sum
$H_n(\{\mu_k\},\{\lambda_k\},{f};z)$. From here by (\ref{ssmm}) we
get the set of the varied moments
\begin{equation}
\label{diff_moments} s_m=m, \quad m\neq n-1,2n-1; \qquad
s_{n-1}=n-1+p,\qquad s_{2n-1}=2n-1+q,
\end{equation}
which are independent of $f$. Consequently,
\begin{equation}
\label{G_n_diff_2_diag}
\hat{G}_n(\lambda):=\sum_{m=0}^n \hat{g}_m\lambda^m=
\left|
\begin{array}{cccccc}
1 & \lambda &  \ldots & \lambda^{n-1} & \lambda^n\\
0 & 1 &  \ldots & n-1+p & n\\
1 & 2 &  \ldots & {n} & {n+1}\\
\ldots &  \ldots & \ldots & \ldots & \ldots\\
n-1+p & n  & \ldots & {2n-2} & {2n-1+q}\\
\end{array}
\right|.
\end{equation}
If for some $p$ and $q$ the generating polynomial
$\hat{G}_n(\lambda)$ has exactly $n$ pairwise distinct roots
$\lambda_1,\ldots,\lambda_n$, then by Theorem~\ref{th1} the varied
interpolation problem becomes regular and
\begin{equation}
\label{REG PROIZ}
zf'(z)=H_n(\{\mu_k\},\{\lambda_k\}, f;z)- p\,
f_{n-1}z^{n-1}-q\,f_{2n-1}z^{2n-1}+O(z^{2n}),
\end{equation}
where $\mu_k$ can be calculated using (\ref{diff_moments}),
(\ref{SRS}) and Lemma~\ref{400}.

\medskip

{\bf \thesection.2. Coefficients of the generating polynomial.} In the case under consideration the
coefficients $\hat{g}_m$ can be written explicitly.
\begin{lemma}
\label{lemma_dif_1} Let $\kappa:=(-1)^{n(n+1)/2}p^{n-3}$. Then for
$n\ge 1$ the coefficients of the polynomial
$(\ref{G_n_diff_2_diag})$ have the form
\begin{align*}
\hat{g}_n&=\kappa p\left(p^2+n(n-1)p +\tfrac{n^2(n^2-1)}{12}\right),\\
\hat{g}_0&=-\kappa \left(p^2q+(2n-1)p^2+(n-1)^2p\,q-\tfrac{n(n^2-1)}{6}p+\tfrac{(n-2)n(n-1)^2}{12}q\right),\quad\quad\quad\;\\
\hat{g}_m&=-\kappa \left((2n-(m+1))p^2-(n-(m+1))p\,q -\tfrac{n(n+1)}{2}\left(\tfrac{n+2}{3}-(m+1)\right)p+\phantom{\tfrac{1}{1}}\right. \\
&\qquad\qquad\qquad\qquad\qquad\left.\phantom{\tfrac{1}{1}}
+\tfrac{n(n-1)}{2}\left(\tfrac{2(n+1)}{3}-(m+1)\right)q\right),\qquad
m=\overline{1,n-1}.
\end{align*}
\end{lemma}
\begin{Proof}
One can verify the identities for $n=1,2$ directly. From now on, $n\ge 3$.

Let us first prove the identity for $\hat{g}_n$ by direct
calculation of the algebraic adjunct  $(-1)^n D$ to $\lambda^n$ in
the determinant (\ref{G_n_diff_2_diag}).

We now show that the characteristic polynomial $P_n(\lambda)=\det
(A-\lambda I)$ of the matrix
$$
A:=\left(
\begin{array}{ccccc}
n-1 & n-2 & n-3 & \ldots & 0 \\
n & n-1 & n-2 & \ldots & 1 \\
\ldots & \ldots & \ldots & \ldots & \ldots \\
2n-2 & 2n-3 & 2n-4 & \ldots & n-1 \\
\end{array}
\right)
$$
has the form
\begin{equation}
\label{charac_poly}
P_n(\lambda)=(-1)^n\lambda^{n-2}\left(\lambda^2-n(n-1)
\lambda+\tfrac{n^2(n^2-1)}{12}\right).
\end{equation}
It is known that for any matrix $B$
\begin{equation}
\label{characht_property_1}
\det (B-\lambda I)=(-1)^n(\lambda^n-b_1\lambda^{n-1}+b_2\lambda^{n-2}+\ldots+(-1)^nb_n),
\end{equation}
where  $b_j$ is the sum of all $j$-rowed diagonal minors of the
matrix $B$ (see, for instance, \cite[\S 3.10]{Jacobson}). In
particular, in terms of the traces of the matrices $B$ and $B^2$ we
have
\begin{equation}
\label{characht_property_2} b_1=\mathrm{Tr}\,B, \qquad
b_2=\tfrac{1}{2}\left((\mathrm{Tr}\,B)^2-\mathrm{Tr}\,B^2\right).
\end{equation}
For our matrix $A$ all minors of the size greater than  two are zero
(as subtracting a row from any other one gives a constant row)
therefore $\rank A=2$. Consequently, the coefficients before the
terms with the powers less than $n-2$ in $\det (A-\lambda I)$ are
zero. Furthermore, it is clear that $\mathrm{Tr}\,A=n(n-1)$. We now
consider the coefficient before $\lambda^{n-2}$.
 It is easily seen (by direct multiplication of the $k$th row by the $k$th column of the matrix $A$) that
$$
\mathrm{Tr}\,A^2=
\sum_{k=1}^n\sum_{m=0}^{n-1}\left((n-1)^2-(m-k+1)^2\right)=
\tfrac{1}{6}n^2(n-1)(5n-7).
$$
It follows that
$$
\tfrac{1}{2}\left((\mathrm{Tr}\,A)^2-\mathrm{Tr}\,A^2\right)
=\tfrac{1}{2}\left(n^2(n-1)^2-\tfrac{1}{6}n^2(n-1)(5n-7)\right)
=\tfrac{n^2(n^2-1)}{12}.
$$
This completes the proof of the formula (\ref{charac_poly}).

Let us return to the determinant $D$. Its matrix is mirror symmetric
with respect to $A$ (i.e. its columns are placed in a reversed
order) and can be obtained by right multiplication  of $A$ by the
anti-diagonal identity matrix. As is known, the determinant of the
$n\times n$ anti-diagonal identity matrix  is equal to
$(-1)^{n(n-1)/2}$. Hence
\begin{equation}
\label{DETER-D} D=(-1)^{\tfrac{n(n-1)}{2}}\det (A+p
I)=(-1)^{\tfrac{n(n-1)}{2}}P_n(-p).
\end{equation}
This and (\ref{charac_poly}) yield the desired formula for $\hat{g}_n=(-1)^n D$, $n\ge 3$.

For convenience, we introduce the set of three elements
\begin{equation}
\label{Pi} \hat{\Pi}:=\left\{0; \tfrac{n}{2}\left(1-n+d_n\right);
\tfrac{n}{2}\left(1-n-d_n\right)\right\},\qquad
d_n:=\sqrt{\tfrac{2}{3}(n-1)(n-3)}.
\end{equation}
Note that if $p\notin \hat{\Pi}$, then $\hat{g}_n\neq 0$. Assume
that $p\notin \hat{\Pi}$ and $\hat{g}_n$ are known. Then we can
determine the desired identities for the other  $n$ coefficients
from the following system of  $n$ linear equations:
\begin{equation}
\label{Newton_formula} \sum_{m=0}^n s_{v-m}\hat g_{n-m}=0,\qquad
v=\overline{n,2n-1}.
\end{equation}
The equation for each $v$ can be obtained by summarizing the
products of the algebraic adjuncts to the elements of the first row
of the determinant (\ref{G_n_diff_2_diag}) (generally, of the
determinant (\ref{G_n})) and the corresponding elements of the
$(v+2)$th row. The linear system (\ref{Newton_formula}) with
unknowns $\hat{g}_0,\ldots,\hat{g}_{n-1}$ has a non-singular matrix
(its determinant is equal to the coefficient $\hat{g}_n\ne 0$,
$p\notin \hat{\Pi}$), hence it has a unique solution.

Consequently, in order to complete the proof of
Lemma~\ref{lemma_dif_1} it is sufficient to verify
(\ref{Newton_formula}) by direct substitution of the moments
(\ref{diff_moments}) and coefficients given in
Lemma~\ref{lemma_dif_1}. This verification is quite simple and can
be reduced to calculation of the sums $\sum_{m=0}^{n}m^\nu$ for
$\nu=1,2$, so we do not dwell on it.

Finally, let $p\in \hat{\Pi}$. The case $p=0$ is not interesting as
then we have $\hat{G}_n(\lambda)\equiv 0$ (see
(\ref{G_n_diff_2_diag}) for $n\ge 3$).  In the case $p\in
\hat{\Pi}\setminus \{0\}$ the system (\ref{Newton_formula}) is
homogeneous and has a singular matrix hence it has infinitely many
non-zero solutions, one of those is given in Lemma~\ref{lemma_dif_1}
(we use that the coefficients $\hat{g}_m$ depend on $p$
continuously).
\end{Proof}

\medskip

{\bf \thesection.3. Factorization of the coefficients of the generating polynomial.} The following lemma about factorization of the
coefficients of $\hat{G}_n$ is fundamental since it enables us to use
Theorem~\ref{th1} in the problem under consideration.
\begin{lemma}
\label{lemma_factor} Let $n\ge 3$, $p\notin \hat{\Pi}$ $($see $(\ref{Pi})$$)$ and
\begin{equation}
\label{PARAM p q} q=q_0(p):=-2\,{\frac {p\left
(3\,p+{n}^{2}-1\right )}{\left ({n}-1 \right )\left (n-2\right
)}}.
\end{equation}
Then the ratios $\hat{g}_m/\hat{g}_n$ for $m=\overline{1,n}$ are
independent of  $p$ and $\hat{g}_0/\hat{g}_n$ depend on $p$
linearly. More precisely, the generating polynomial has the form
 \begin{equation}
  \label{FACTOR G}
\hat{G}_n(\lambda)=\hat{g}_n \left({\lambda}^{n}-
\frac{6\lambda\left(\lambda^{n-1}-(n-1)\lambda+n-2\right)}{\left (n-
1\right)\left(n-2\right)\left (\lambda-1\right)^{2}} +2+{\frac
{6\,p}{\left (n-1\right )\left (n-2\right)}} \right).
 \end{equation}
There exists an arbitrarily small variation of the parameter $p$
such that all the roots of the polynomial $\hat{G}_n$ are pairwise
distinct.
\end{lemma}
\begin{Proof} Taking into account Lemma~\ref{lemma_dif_1}, if we solve the equation $\hat{g}_{n-1}=0$, being linear with respect to
$q$, then we get (\ref{PARAM p q}). Substitution of the expression
for $q$ into the other coefficients gives
$$
\hat{g}_0=\hat{g}_n\left(2+{\frac {6\,p}{\left (n-1\right )\left (n-2\right
)}}\right),\qquad \hat{g}_m=-6\,\hat{g}_n{\frac {n-1-m}{\left (n-1\right
)\left ({n}-2 \right )}},\qquad m=\overline{1,n-1},
$$
where $\hat{g}_n=(-1)^{\tfrac{n(n+1)}{2}}p^{n-2}\left(p^2+n(n-1)p
+\tfrac{n^2(n^2-1)}{12}\right)\neq 0$ as $p\notin \Pi$, therefore
$$
\hat{G}_n(\lambda)=\hat{g}_n\, \left ({\lambda}^{n}-\frac{6}{\left (n-
1\right )\left (n-2\right )}\sum _{m=1}^{n-1}\,{\left (n-m-1\right
){\lambda}^{m}}+2+{\frac {6\,p}{\left (n-1\right )\left (n-2\right
)}} \right ),
$$
which  yields (\ref{FACTOR G}) after calculation of the sum.
The conclusion about the simplicity of the roots
 follows immediately from Lemma~\ref{RAZDEL}.
\end{Proof}

\medskip

{\bf \thesection.4. Main theorem about numerical differentiation by amplitude and frequency operators.} From (\ref{corollary_form_Newton++++}) we get
\begin{equation}
\label{corollary_form_Newton} S_{2n}
=-\frac{1}{\hat{g}_n}\sum_{m=0}^{n-1} s_{n+m}\hat{g}_{m}.
\end{equation}
After substituting the moments (\ref{diff_moments}) and the
coefficients from Lemma~\ref{lemma_factor} into
(\ref{corollary_form_Newton}), direct calculation gives the
expression for the $2n$th moment:
$$
S_{2n} = 2n-C_n(p),\qquad C_n(p):=\frac{6np}{(n-1)(n-2)}.
$$

Therefore the remainder $r_n$, which was denoted in (\ref{REG
PROIZ}) just as $O(z^{2n})$, has the following more precise form:
\begin{equation}
\label{r_n_diff} \tilde{f}(z)-\sum_{k=1}^n\mu_k f(\lambda_k z)=
\sum_{m=2n}^{\infty}(m-S_m)f_mz^m= C_n(p)f_{2n}z^{2n}+O(z^{2n+1}).
\end{equation}

From the foregoing, we get the following statement.

\begin{theorem}
\label{th3}  Given $n\ge 3$, any $p_0\in \mathbb C$ and arbitrarily
small $\varepsilon>0$, there exists a value of the parameter $p$,
$|p-p_0|\le \varepsilon$, $p\notin \hat{\Pi}$ $($see $(\ref{Pi}))$,
such that
\begin{equation}
\label{diff_formula_common} zf'(z)=\sum_{k=1}^n\mu_k f(\lambda_k
z)-pf_{n-1}z^{n-1}- qf_{2n-1}z^{2n-1}+r_n(z),
\end{equation}
where $r_n(z)=O(z^{2n})$ is the form $(\ref{r_n_diff})$ and
$q=q_0(p)$ $($see $(\ref{PARAM p q})$$)$. Moreover, the frequencies
$\lambda_k$ are the pairwise distinct roots of the polynomial
$(\ref{FACTOR G})$ and the amplitudes $\mu_k$ are determined
uniquely by Lemma~\ref{400}. Furthermore, $\mu_k=\mu_k(p,n)$ and
$\lambda_k=\lambda_k(p,n)$, so they are independent of the function
$f$ and universal in this sense.

The interpolation formula is exact for the polynomials of degree
$\le 2n-1$, i.e. ${r_n(z)\equiv 0}$ in $(\ref{diff_formula_common})$
for the polynomials $f$ such that $\deg f\le 2n-1$.
\end{theorem}

Note that the method, which we consider in this section, can be
easily extended to interpolation of the functions $z^\nu
f^{(\nu)}(z)$, $\nu\ge 2$, hence one can obtain the formulas for
numerical differentiation of higher order.

\medskip
 \textbf{\thesection.5. Remarks and examples.}
 We now make several remarks about practical applications of Theorem~\ref{th3}.
  The remainder in the formula (\ref{diff_formula_common})
  is of  quite high infinitesimal order, $O(z^{2n})$, and this is achieved by knowing only $n$ values of $f$
 and two fixed values of its derivatives at $z=0$. Traditional interpolation approaches with such a number of
 known values usually have remainders of order $O(z^{n+2})$. In other words, the formula
$(\ref{diff_formula_common} )$ is exact for the polynomials of
degree $\le 2n-1$, whereas usual $(n+2)$-point interpolation
formulas are exact only for the polynomials of degree $\le n+1$.

Another important feature of the formula (\ref{diff_formula_common})
is that the variable interpolation nodes $\lambda_k z$ depend only
on the point $z$, where we calculate $zf'(z)$, and are independent
of $f$ (in this sense the amplitudes $\mu_k=\mu_k(p,n)$ and
frequencies $\lambda_k=\lambda_k(p,n)$ are universal for the whole
class of analytic functions).

It is seen from the formula (\ref{diff_formula_common}) that its
precision strongly depends on the precision of the values
$f_{\nu}=f^{(\nu)}(0)/\nu!$, $\nu=n-1,2n-1$ (of course we assume
that the values of the function $f$ are known).
 In several particular cases this difficulty can be overcome. For instance, if it is known a priori that
the function $f$ is even (odd), then there is no necessity in
calculation of $f_{n-1}$ and $f_{2n-1}$ for even (odd) $n$ since
then $pf_{n-1}z^{n-1}+ qf_{2n-1}z^{2n-1}\equiv 0$
($pf_{n-1}z^{n-1}\equiv 0$) and the local precision is $O(z^{2n})$
($O(z^{2n-1})$). In the more general case, when $f$ is even or odd,
the formula (\ref{diff_formula_common}) can be applied to the even
auxiliary function $\omega(z)=f^2(z)$ for even $n$. Then the
corresponding coefficients $\omega_{n-1}=\omega_{2n-1}\equiv 0$ and
$$
 2zf(z)f'(z)=\sum_{k=1}^n\mu_k
f^2(\lambda_k z)+O(z^{2n}).
$$
In the most general case, if we want to use (\ref{diff_formula_common}) systematically, it is necessary to calculate
the regularizing binomial $pf_{n-1}z^{n-1}+qf_{2n-1}z^{2n-1}$  for
each {\it fixed} function~$f$. For this purpose some known formulas
for numerical differentiation of analytic functions at $z=0$ can be
used. For example, one can use several high-accuracy formulas for
$f_{\nu}=f^{(\nu)}(0)/\nu!$, obtained in \cite{Schmeisser}.

Note that, in contrast to the formulas for calculation of Taylor
coefficients as, for instance, in \cite{Schmeisser}, the formula
 (\ref{diff_formula_common}) works well only in a deleted neighbourhood of the point $z=0$.
Below we cite several known interpolation formulas for numerical
differentiation, being close in form to (\ref{diff_formula_common}).

In \cite{Ash_Janson_Jones} the following
$n$-point formulas for numerical differentiation of real functions
were obtained:
\begin{equation}
\label{1} f_\nu x^{\nu}= \sum_{k=1}^n\mu_k f(\lambda_k
x)+O(x^n),\qquad {\nu}=1,2,\qquad n\ge {\nu}+1,
\end{equation}
where $\lambda_k x$, $|\lambda_k-\lambda_j|>1$, are real nodes,
minimizing the generalized power sums
$S_v=\sum_{k=1}^n\mu_k\lambda_k^v$ for $v\ge n+1$ (this corresponds to
minimization of the remainder). In \cite{Salzer} interpolation
formulas of the Lagrange type for numerical differentiation were
constructed on basis of special non-uniformly distributed nodes,
also minimizing the remainder.

Formulas of the form (\ref{1}) for analytic functions  were obtained
in \cite{Lyness} via contour integrals. Moreover, it was shown there
that for $\lambda_k=\exp (2\pi i (k-1)/n)$ and appropriate $\mu_k$
their formulas were exact for the polynomials of degree $\le
n+{\nu}-1$. This result can be extended. Indeed, by direct substitution
 (see the discussion near the formulas (\ref{alpha-beta}) and (\ref{TAU}))
one can check that for any non-vanishing parameters $p$ and $q$ and
any integer $0\le \nu\le n-1$ we have
 $$
pf_\nu z^\nu+qf_{\nu+n}z^{\nu+n}=\sum_{k=1}^n\mu_k f\left(\lambda_k
z\right)+ O(z^{2n+\nu}),\quad
\lambda_k=\left(\frac{q}{p}\right)^{1/n}\exp\left(\frac{2\pi
(k-1)i}{n}\right),
$$
where $\mu_k=(-1)^{n-\nu-1}\lambda_k \sigma_{n-\nu-1}^{(k)}\,p^2
\,(qn)^{-1}$ $($see Lemma~$\ref{400})$ and in $\lambda_k$ one can
take any value of the root.

The following formula for analytic functions $f$ is contained in
\cite{Dan2008}:
$$
f_\nu z^{\nu}= \sum_{k=1}^{(\nu+1)
N}\lambda_kf(\lambda_kz)+O(z^{n}), \qquad
N=\left[\tfrac{n}{\nu+1}\right], \qquad n>\nu+6.
$$
Here the numbers $\lambda_k$ do not depend on $f$ and are non-zero
roots of the polynomial $P_n(\lambda)=\sum_{k=0}^N (-1)^k
\lambda^{n-k-\nu k}/((\nu+1)^k k!)$ (see estimates for the remainder
 in \cite{Dan2008}). Other results of this type were also obtained in
\cite{Chu2010,Fryantsev}.
\begin{example}[{\bf numerical differentiation}]
Let $n=4$ and $p=-1$. Then by (\ref{PARAM p q}) and (\ref{FACTOR G})
we get $q=4$ and $\hat{G}_4(\lambda)=
9\left(\lambda^4-\lambda^2-2\,\lambda+1\right)$. The roots of the
generating polynomial $\hat{G}_4$ are
$$
\lambda_1\approx 1.38647,\quad \lambda_2\approx0.42578,\quad
\lambda_{3,4}\approx-0.90612\pm 0.93427\,i.
$$
Lemma~\ref{400} gives
$$
 \mu_1\approx0.967276,\quad \mu_2\approx-0.79945,\quad
\mu_{3,4}\approx-0.08390\pm 0.08175\, i,
$$
and the formula (\ref{diff_formula_common}) takes the form
\begin{equation}
\label{diff_example_n=4} zf'(z)=\sum_{k=1}^4\mu_k f(\lambda_k
z)+f_{3}z^{3}-4f_{7}z^{7}+r_4(z),\quad r_4(z)=-4 f_8 z^8-9 f_9
z^9+\cdots.
\end{equation}

For instance, set $f(z)=(z+2)^{-1}$ (then $f_3=-1/16$,
$f_7=-1/256$). Calculations in Maple show that the error~of
~(\ref{diff_example_n=4}) does not exceed $10^{-4}$ for $z\in
[-0.5,0.5]$. For $n=7$ and $n=10$ the corresponding errors are less
than $10^{-8}$ and $10^{-12}$, correspondingly.

Now consider the Bessel function  $f=J_0$ from (\ref{J_0000}). This
function is even and consequently the regularizing binomial
$pf_{n-1}z^{n-1}+ qf_{2n-1}z^{2n-1}$ vanishes for even $n$.
Therefore from (\ref{diff_example_n=4}) we get
$$
zJ'_0(z)\approx \sum_{k=1}^4 \mu_kJ_0(\lambda_kz).
$$
The error of the approximant does not exceed  $10^{-4}$ for
$z\in[-1,1]$. For $n=6$ and $n=8$ the corresponding errors are less than $10^{-9}$ and $10^{-14}$.
\end{example}

\textbf{\thesection.6. Some estimates.} Absolute values of the amplitudes $\mu_k$
and frequencies $\lambda_k$ play an important role in calculations
by the formula (\ref{diff_formula_common}). We now estimate the
frequencies.
\begin{lemma}
\label{500} For the roots $\lambda_k$ of the polynomial $(\ref{FACTOR G})$ we have
$$
|\lambda_k|\le 1+\frac{O(1)}{\sqrt{n}}, \qquad O(1)>0, \qquad n\to \infty.
$$
More precisely, given $n\ge 3$,
$$
|\lambda_k|\le
\Lambda:=\left(2\delta\right)^{\frac{3}{\sqrt{n-2}}},\qquad
\delta:=1+\frac{3|p|}{(n-1)(n-2)}, \qquad p\notin \Pi.
$$
\end{lemma}
\begin{Proof}
First, we estimate the absolute value of the sum of the last three
terms in the brackets in $(\ref{FACTOR G})$:
\begin{align*}
V:&=\left|-\frac{6\lambda\left({\lambda}^{n-1}-(n-1){\lambda}+n-2\right)}{\left(n-1\right)\left(n-2\right)\left (\lambda-1\right)^{2}}
+2+{\frac {6\,p}{\left (n-1\right )\left
(n-2\right )}}\right|\\
&\le |\lambda|^n \left(\frac{6\left(1+(n-1)|\lambda|^{2-n}+(n-2)|\lambda|^{1-n}\right)}{(n-1)(n-2)\left(|\lambda|-1\right)^{2}}
 +\frac{2\delta}{|\lambda|^n}\right).
\end{align*}
It is easily seen that
$\left(2\delta\right)^{\frac{3}{\sqrt{n-2}}}-1\ge \tfrac{3\ln
2}{\sqrt{n-2}}$ for $\delta>1$. Therefore substituting
$|\lambda|=\Lambda$ into the latter expression yields
$$
V\le |\lambda|^n \left(\frac{6
\cdot\left(1+(2n-3)/(2\delta)^{3\sqrt{n-2}}\right)}{(3\ln 2)^2
\cdot(n-1)}
 +\frac{1}{(2\delta)^{\frac{3n}{\sqrt{n-2}}-1}}\right).
$$
It is also easy to check that for $n\ge 3$ and $\delta>1$ the
expression in the brackets is less than one, so $V<\Lambda^n$ for
the above-mentioned $|\lambda|=\Lambda$. This and Rouch\'{e}'s
theorem imply that all the roots of the polynomial (\ref{FACTOR G})
lie in the disc $|\lambda|< \Lambda$.
\end{Proof}

One can use Lemmas~\ref{400} and \ref{500} to obtain estimates for
the amplitudes $\mu_k$ but this problem needs more delicate
 analysis and we do not dwell on it in the present paper.

\section{Numerical extrapolation by amplitude and frequency operators}
\label{par-extrap}
 {\bf \thesection.1. Statement of the problem.} Let us briefly describe the
idea of the extrapolation. Let $a>0$, $p,q\in \mathbb{R}$ and $f$ be
a function analytic in a disc $|z|<\rho$, $\rho>0$. Consider the
 problem of multiple interpolation of the function $f(az)$ in a
neighbourhood of $z=0$  by the amplitude and frequency operator
${H}_n(\{\mu_k\},\{\lambda_k\},f;z)$, where $f$ is chosen as a basis
function. As in the case of differentiation, we get a non-regular
discrete moment problem with $s_m=a^m$. To regularize it, we
introduce the varied function
\begin{equation}
\label{ZADACHA extrapol}
\tilde{f}(z):=f(az)+pf_{n-1}z^{n-1}+qf_{2n-1}z^{2n-1}
\end{equation}
with some parameters $p$ and $q$, being non-zero simultaneously. By
the same approach that we used at the beginning of
Section~\ref{par-diff}, in order to construct the interpolating sum
${H}_n(\{\mu_k\},\{\lambda_k\},f;z)$, we find the sequence of varied
moments
\begin{equation}
\label{moments_extropal} s_k=a^{k}, \;\; k\neq n-1,2n-1, \;\;
k\in\mathbb{N}_0; \quad s_{n-1}=a^{n-1}+p,\quad s_{2n-1}=a^{2n-1}+q,
\end{equation}
and construct the generating polynomial
\begin{equation}
\label{G_n_extrapol_p_q} \check{G}_n(\lambda):=\sum_{m=0}^n
\check{g}_m\lambda^m= \left|
\begin{array}{cccccc}
1 & \lambda &  \ldots & \lambda^{n-1} & \lambda^n\\
1 & a &  \ldots & a^{n-1}+p & a^n\\
a & a^2 &  \ldots & a^{n} & a^{n+1}\\
\cdots &  \cdots & \cdots & \cdots & \cdots\\
a^{n-1}+p & a^n  & \ldots & a^{2n-2} & a^{2n-1}+q\\
\end{array}
\right|.
\end{equation}
If for some $a> 0$, $p$ and $q$ the polynomial $\check{G}_n$ is of
degree $n$ and all its  roots $\lambda_1,\ldots,\lambda_n$ are
pairwise distinct, then by Theorem~\ref{th1} the varied problem for
the function
 (\ref{ZADACHA extrapol}) is regularly solvable, so the following interpolation formula
 holds:
\begin{equation}
\label{100} f(az)=H_n(\{\mu_k\},\{\lambda_k\}, f;z)- p\,
f_{n-1}z^{n-1}-q\,f_{2n-1}z^{2n-1}+O(z^{2n}).
\end{equation}
(Of course, we assume that all the arguments of the function $f$ lie
in the disc ${|z|<\rho}$, where it is analytic.) Suppose also that
the inequalities  ${|\lambda_k|<\delta a}$ with some ${\delta\in
(0,1)}$ are valid for all $k=\overline{1,n}$. Then it is natural to
call the formula (\ref{100}) {\it extrapolational} as the values of
the function $f$ at the points $\zeta=az$ are approximated by the
values of this function at the points $\lambda_kz$, belonging to the
disc $\{\xi: |\xi|<\delta |\zeta|\}$, $\delta <1$. In the present
section we will obtain such an extrapolation formula and a quantitative
estimate for its remainder.

We start with a formal description. As in Section~\ref{par-diff}, we
first analyse the coefficients and roots of the
polynomial~$\check{G}_n$ of the form (\ref{G_n_extrapol_p_q}).

\medskip

{\bf \thesection.2. Coefficients of the generating polynomial.} The following statement gives an explicit form
of the coefficients~$\check{g}_m$.
\begin{lemma}
\label{lemma_extrapol_1} Let $\kappa:=(-1)^{n(n+1)/2}p^{n-2}$,
$p\neq 0,\;-na^{n-1}$. The polynomial $\check{G}_n$ has the
following coefficients:
\begin{equation}
\label{koef_extrapol}
\begin{array}{c}
\check{g}_n=\kappa p\left(na^{n-1}+p\right),\qquad \check{g}_0=-\kappa\left(a^{2n-1}p+(n-1)a^{n-1}q+p\,q\right), \\
\check{g}_m=-\kappa a^{n-1-m}\left(a^np-q\right),\qquad
m=\overline{1,n-1}.
\end{array}
\end{equation}
\end{lemma}
\begin{Proof}
The method of proof is the same as in Lemma~\ref{lemma_dif_1}. We
first prove the identity for ${\check{g}}_n$ by direct computation
of the algebraic adjunct $(-1)^n D$ to the element $\lambda^n$ in
the determinant (\ref{G_n_extrapol_p_q}). For this we show that
given the matrix
$$
A:=\left(
\begin{array}{ccccc}
a^{n-1} & a^{n-2} & a^{n-3} & \ldots & 1 \\
a^n & a^{n-1} & a^{n-2} & \ldots & a \\
\ldots & \ldots & \ldots & \ldots & \ldots \\
a^{2n-2} & a^{2n-3} & a^{2n-4} & \ldots & a^{n-1} \\
\end{array}
\right),
$$
the characteristic polynomial $P_n(\lambda)=\det (A-\lambda I)$ has the form
\begin{equation}
\label{charac_poly_2}
P_n(\lambda)=(-1)^n\lambda^{n-1}\left(\lambda-na^{n-1}\right).
\end{equation}
Indeed, in this case $\rank A=1$ as any two rows of   $A$ are
proportional. Therefore, by (\ref{characht_property_1}) and
(\ref{characht_property_2}) the coefficients before the terms with
the powers less than  ${n-1}$ are zero in $\det (A-\lambda I)$. To
prove (\ref{charac_poly_2}), it remains to notice that $\mathrm{Tr}
A=na^{n-1}$.

Now we return to the determinant $D$. In the same way as in Lemma~\ref{lemma_dif_1},
 (\ref{DETER-D}) and   (\ref{charac_poly_2}) yield the desired formula for ${\check{g}}_n$.

If we suppose that ${\check{g}}_n$ are known, then the other $n$
coefficients as in the above-mentioned case of differentiation can
be found from the system (\ref{Newton_formula}) of $n$ linear
equations (with the exchange of $\hat{\phantom{o}}$ by
$\check{\phantom{o}}$). This system with respect to unknowns
${\check{g}}_0,\ldots,{\check{g}}_{n-1}$ has a unique solution for
$p\neq 0,\;-na^{n-1}$ as its determinant is equal to
${\check{g}}_n\ne 0$. Thus it suffices to check
(\ref{Newton_formula}) by direct substitution of the values
(\ref{moments_extropal}) and (\ref{koef_extrapol}). This is quite
easy and thus we do not dwell on it.
\end{Proof}

\medskip

{\bf \thesection.3. Roots of the generating polynomial.}  From now on we suppose that $q=0$; then, as we
will see below, $\check{G}_n$ has exactly $n$ pairwise distinct
roots and its coefficients can be calculated by quite simple
formulas.
\begin{lemma}
Let $p>0$, $q=0$ and $a>0$. Then
\begin{equation}
 \label{COEFF G_n_extrapol_p_q_special}
{\check{g}}_n=\kappa p\left(na^{n-1}+p\right)\neq 0,\qquad
{\check{g}}_m=-\kappa p \,a^{2n-1-m},\qquad m=\overline{0,n-1},
\end{equation} and the generating polynomial has the form
\begin{equation}
 \label{G_n_extrapol_p_q_special}
\check{G}_n(\lambda)=
{\check{g}}_n\left(\lambda^n-\frac{a^{2n-1}}{na^{n-1}+p}\sum_{m=0}^{n-1}
\frac{\lambda^m}{a^m}\right)={\check{g}}_n
\left(\lambda^n-\frac{a^{n}}{na^{n-1}+p}
\frac{\lambda^n-a^n}{\lambda-a}\right).
\end{equation}
Moreover, $\check{G}_n$ has exactly $n$  pairwise distinct roots.
\end{lemma}
\begin{Proof}
The representation (\ref{G_n_extrapol_p_q_special}) can be obtained
by direct substituting $q=0$ into the formulas
(\ref{koef_extrapol}) from Lemma~\ref{lemma_extrapol_1} for $p>0$.
It remains to show that the polynomial
(\ref{G_n_extrapol_p_q_special}) has no multiple roots. We rewrite
$\check{G}_n$ in the form
$$
\check{G}_n(\lambda)={\check{g}}_n
\frac{P_{n+1}(\lambda)}{\lambda-a},\qquad
P_{n+1}(\lambda):=\lambda^{n+1}-a\
\left(1+\frac{a^{n-1}}{na^{n-1}+p}\right)\lambda^n+\frac{a^{2n}}{na^{n-1}+p}.
$$
The set of roots of the polynomial $P_{n+1}$ contains all the roots
of the polynomial $\check{G}_n$ and one more root $\lambda=a$. If
the polynomial $P_{n+1}$ had a multiple root, then this root would be also
a root of its derivative  $P'_{n+1}$. However,
$$
P'_{n+1}(\lambda)=(n+1)\lambda^{n-1}\left(\lambda-\lambda^{*}\right),\qquad
\lambda^{*}:=\frac{a
n}{n+1}\left(1+\frac{a^{n-1}}{na^{n-1}+p}\right),
$$
and, as it can be easily seen, in both roots of the derivative
$P'_{n+1}$, namely, $0$ and $\lambda^{*}$, the polynomial $P_{n+1}$
does not vanish (for $p>0$), more precisely, $P_{n+1}(0)>0$,
$P_{n+1}(\lambda^{*})<0$. Consequently, both $P_{n+1}$ and $\check{G}_n$  have no multiple roots.
\end{Proof}

\medskip

{\bf \thesection.4. Estimates of the roots of the generating polynomial.}  We now estimate the roots of the polynomial
(\ref{G_n_extrapol_p_q_special}) assuming $q=0$ as above.
\begin{lemma}
\label{lemma_extrapol_2} Given $p>0$ and $a>0$, for the roots
$\lambda_k$ of the polynomial $\check{G}_n$ we have
\begin{equation}
 \label{COR-G_n}
|\lambda_k|<\delta a<a,\qquad
\delta=\delta(n,a,p):=\left(1+\frac{p}{na^{n-1}}\right)^{-1/n},
\qquad k=\overline{1,n}.
\end{equation}
\end{lemma}
\begin{Proof} For $|\lambda|=\delta a$ we get
$$
\frac{a^{2n-1}}{na^{n-1}+p}\sum_{m=0}^{n-1}\left|\frac{\lambda}{a}\right|^m=
\frac{a^{2n-1}}{na^{n-1}+p}\sum_{m=0}^{n-1}\left|\delta\right|^m<\frac{na^{n-1}}{na^{n-1}+p}\,a^n=
(\delta a)^n=|\lambda|^n.
$$
From this, taking into account the first identity in
$(\ref{G_n_extrapol_p_q_special})$ and Rouch\'{e}'s theorem, we
conclude that all $n$ roots of the polynomial $\check{G}_n$ belong
to the disc $|\lambda|<\delta a$.
\end{Proof}
\begin{remark}
\label{primech_ext} We now mention several properties of
$\delta=\delta(n,a,p)$, which plays a key role in the process of
extrapolation. For a fixed $n$ we obviously have
\begin{equation*}
 \label{COR-G_n-asymp}
\delta\in (0,1),\qquad
\delta=\left(\frac{n}{p}\right)^{1/n}a^{1-\tfrac{1}{n}}
\left(1-O\left(\frac{a^{n-1}}{p}\right)\right),\qquad
\frac{a^{n-1}}{p}\to 0,
 \end{equation*}
where $O\left(a^{n-1}/p\right)$ is a positive real value.

From this it follows, in particular, that if the fraction
$a^{n-1}/p$ decreases, then all the roots $\lambda_k$ come closer to
the origin. For example, for a fixed $p$ and $a\to 0$ the largest
absolute value of $\lambda_k$ is bounded by a value of order
$a^{2-1/n}$.
\end{remark}

{\bf \thesection.5. Main theorem about numerical extrapolation by amplitude and frequency operators. Remarks and examples.}  We now aim to estimate the extrapolation
remainder (\ref{100}) and then formulate the main result of the
section. To do so, we need estimates for the generalized power sums
$S_v$, $v\ge 2n$, taking into account the sums  (\ref{moments_extropal}) with indexes $\le 2n-1$ and $q=0$ as before.
\begin{lemma}
\label{lemma_extrapol_3} For $p>0$ and $a>0$ the following
inequalities hold:
\begin{equation}
 \label{200}
0\le S_v\le a^v, \qquad v\ge 2n.
 \end{equation}
\end{lemma}
\begin{Proof}
We prove this by induction, based on (\ref{COEFF
G_n_extrapol_p_q_special}) and the identities
~(\ref{corollary_form_Newton++++}) for $v\ge 2n$. For $v=2n$ from
(\ref{moments_extropal}) and (\ref{COEFF G_n_extrapol_p_q_special})
we get
$$
S_{2n} = \frac{na^{3n-1}}{na^{n-1}+p}\le a^{2n}.
$$
Furthermore, suppose that the inequality $S_v\le a^v$ is valid for
all $v=\overline{2n,N}$ (hence for all $v=\overline{0,N}$). Under
this assumption, we obtain
$$
S_{N+1} = \frac{a^{2n-1}}{na^{n-1}+p}\sum_{m=0}^{n-1}
\frac{S_{N-n+m+1}}{a^m}\le \frac{na^{N+n}}{na^{n-1}+p}\le a^{N+1},
$$
which completes the proof.
\end{Proof}

Now, using (\ref{200}), we can estimate the extrapolation remainder
(\ref{100}), where $q=0$:
$$
|r_n(z)|=
\left|\sum_{m=2n}^{\infty}(a^m-S_m)f_mz^m\right|
\le\sum_{m=2n}^{\infty}|f_m||az|^m.
$$
Note that this estimate is independent of $p$.

\medskip

 Summarizing, we formulate the main result of this
section.
\begin{theorem}
\label{th4} Let $f$ be analytic in the disc  $|z|<\rho$, $a>0$,
$p>0$. Then for $|z|<\rho/a$ the following extrapolation formula
holds:
\begin{equation}
\label{extrapol_formula} f(az)=\sum_{k=1}^n\mu_k f(\lambda_k
z)-pf_{n-1}z^{n-1}+r_n(z),\qquad |r_n(z)|\le
\sum_{m=2n}^{\infty}|f_m||az|^m,
\end{equation}
where the frequencies $\lambda_k$ are the pairwise distinct roots of
the polynomial $(\ref{G_n_extrapol_p_q_special})$, and
$($see~$(\ref{COR-G_n}))$
$$
|\lambda_k z|<\delta a|z|<a|z|,\qquad
\delta=\left(1+\frac{p}{na^{n-1}}\right)^{-1/n}.
$$
The amplitudes $\mu_k$ are uniquely determined by Lemma~\ref{400}.
Moreover, $\lambda_k=\lambda_k(a,p,n)$ and $\mu_k=\mu_k(a,p,n)$, so they are independent of the function  $f$ and universal in this sense.

The extrapolation formula $(\ref{extrapol_formula})$ is exact for
the polynomials of degree $2n-1$, i.e. $r_n(z)\equiv 0$ for the
polynomials $f$ such that $\deg f\le 2n-1$.
\end{theorem}
\begin{remark}
If $f_{n-1}=0$, then the extrapolation formula has a particular
simple form and high degree of accuracy. For instance, for even
(odd) functions and even (odd) $n$ the identity from
(\ref{extrapol_formula}) has the form
$$
f(az)=\sum_{k=1}^n\mu_k f(\lambda_k z)+r_n(z).
$$

Calculation of $f_{n-1}$ in the general case was
discussed at the end of Section~\ref{par-diff}.
\end{remark}
\begin{remark}
We emphasize that the extrapolation character of the formula
(\ref{extrapol_formula}) is specified by a proper choice of the
parameters $p$ and $q$ in the problem (\ref{moments_extropal}). For
instance, for $p=q=0$ this problem is also solvable (but
non-regular): one can take $\mu_1=1$ and $\lambda_1=a$ with the rest
of parameters $\mu_k$, being equal to zero, and any $\lambda_k$.
However, in this case (\ref{extrapol_formula}) becomes a
trivial identity.

If one takes, for example, $p>0$ and $q=a^n p-\check{g}_n/\kappa$, then the problem
 (\ref{moments_extropal}) turns out to be regular and, moreover, the
 coefficients of the generating polynomial $\check{G}_n$ can be also calculated easily:
$$
\check{g}_0=-\check{g}_np_0, \qquad
p_0:=\left(na^{n-1}+p\right)\left(\kappa p
\,a^n/\check{g}_n-1\right)+a^{n-1},
$$
$$
\check{g}_m=-a^{n-1-m}\check{g}_n,\qquad m=\overline{0,n-1}.
$$
However, in this case  there is no extrapolation since the inequalities  $|\lambda_k|<a$ are not valid anymore and we just have an interpolation formula of the form (\ref{extrapol_formula}).
\end{remark}
\begin{remark}
Note that for different $a$ one has different extrapolation formulas
of the form (\ref{extrapol_formula}), not arising from each other.
In particular, they cannot be reduced to the case $a=1$ by linear
substitution. Indeed, substituting $t=az$ into
(\ref{extrapol_formula}) gives
\begin{equation}
\label{extrapol_formula_substitute} f(t)=\sum_{k=1}^n\mu_k(a)
f\left(\tilde\lambda_k(a)
t\right)-pf_{n-1}(t/a)^{n-1}+r_n(t/a),\qquad
\tilde\lambda_k(a)=\lambda_k(a)/a,
\end{equation}
where by Lemma~\ref{lemma_extrapol_2} for any $p>0$ and $a>0$
$$
|\tilde\lambda_k(a)|\le
\delta(n,a,p)=\left(1+\frac{p}{na^{n-1}}\right)^{-1/n}<1,
$$
i.e. $a$ does not disappear and still is a controlling parameter. If
$a<1$, then for any fixed $p>0$ we get
$$
|\tilde{\lambda}_k(a)|<\delta(n,a,p)\to a,\qquad n\to\infty,
$$
Thus, for large $n$ all the arguments $\tilde{\lambda}_k(a)\,t$ lie
almost $a$-times closer to the origin than~$t$ (see also
Remark~\ref{primech_ext}).
\end{remark}
\begin{remark}
Realizing the $n$-point simple or multiple extrapolation
(interpolation) on basis of the Lagrange polynomials or other
similar approaches, one usually obtains extrapolation
(interpolation) formulas, being exact for the polynomials of degree
$n-1$ (see, for instance, \cite{Salzer2,DanChu2011,Chu2012}).
However, our extrapolation formula is exact for the polynomials of
degree $\le 2n-1$. It is interesting that the doubling of precision
is gained by adding just one regularizing power term
$pf_{n-1}z^{n-1}$.

We also emphasize that,  due to Remark~\ref{primech_ext}, if
$p\to\infty$ and all other parameters are fixed, then  the
extrapolation nodes tend to the point $z=0$, but at the same time
the theoretical error of extrapolation does not increase (see
(\ref{extrapol_formula})) as is independent of $p$. The same
phenomenon of the convergence of nodes to the origin was noticed in
similar extrapolation problems in \cite{DanChu2011,Chu2012}.
\end{remark}
\begin{example}[{\bf numerical extrapolation}]
Let $n=2$, $a=1/2$ and $p =2$ ($q=0$ as above). Then the generating
polynomial (\ref{G_n_extrapol_p_q_special}) has the form
$$
\check{G}_2(\lambda)=-6\,{{\lambda}}^{2}+\frac{1}{2}\,{\lambda}+\frac{1}{4}.
$$
We find its (pairwise distinct) roots and then determine the
amplitudes by Lemma~\ref{400}:
\begin{equation*}
\label{example1-extrapolation-lambda+} \lambda_{1}=- \frac{1}{6}
\qquad \lambda_{2}=\frac{1}{4},\qquad \mu_{1}=-\frac{27}{5},\qquad
\mu_{2}=\frac{32}{5}.
\end{equation*}
Thus we get the following extrapolation formula (written in the form
(\ref{extrapol_formula_substitute})):
\begin{equation*}
\label{example1-extrapolation++} f(z)\approx -\frac{27}{5}
f\left(-\frac{1}{3}z\right)+\frac{32}{5}f\left(\frac{1}{2}z\right)-4f_1z.
\end{equation*}
For example, for $f(z)=e^z$ the absolute error of this formula does
not exceed $0.002$ in the real segment $z\in [-0.5,0.5]$.

For $n=4$ and $n=8$ and the same parameters $a=1/2$ and $p=2$ the
error of the extrapolation formula  (\ref{extrapol_formula}) for
$e^z$ in $[-1,1]$ does not exceed $10^{-7}$ and $10^{-18}$
correspondingly. Moreover, in both cases the moduli of extrapolation
nodes are less than $0.58|z|$.
\end{example}

\section{On the numerical method for constructing amplitude and frequency operators}
\label{Section6}

As we have already mentioned in Remark~\ref{remark2}, some authors
studied numerical methods for solving the systems~(\ref{SRS_M}). In
\cite{YF} one of such approaches, the method of small residuals in
 overdetermined moment systems, was used for approximation by amplitude and frequency sums~(\ref{gH}) in a neighbourhood of the
point~${z=0}$. Some discussions from this paper (see Theorems~3 and
5 and Remark~1 there) raise the following important question: can
one use the method of small residuals in the context of the Pad\'{e}
interpolation (\ref{2n-interpolation}) and approximation by
amplitude and frequency operators? From the undermentioned
observations it is seen that this method can work well only for the
quite narrow class of consistent systems, but in the general case one
has to give a negative answer to the question. For the discussion
we can consider only $M=2n$ as the case $M>2n$ follows from the
below-mentioned counterexamples by adding equations with
arbitrary right hand sides.

Following \cite{YF}, instead of the system~(\ref{SRS}) we consider
the one with small residuals $\delta:=\{\delta_m\}_{m=0}^{2n-1}$:
\begin{equation}
\label{NVZ} s_m-\sum_{k=1}^n\mu_k\lambda_k^m=\delta_m,\qquad
m=\overline{0,2n-1}, \qquad |\delta|:=\max_m|\delta_m|\le
\varepsilon.
\end{equation}
It is not difficult to show that for an arbitrarily small
$\varepsilon>0$ one can choose the residuals $\delta$, $|\delta|\le
\varepsilon$, such that the system (\ref{NVZ}) is solvable (both for
 consistent and inconsistent systems (\ref{SRS})). Furthermore, using
this solution, one can construct the corresponding amplitude and
frequency sum $H_n(\delta;z)$ of the form (\ref{gH}) such that
\begin{equation}
\label{NVZ+} f(z)-H_n(\delta;z)=\sum_{m=0}^{2n-1} h_m\delta_m
z^m+B_n(\delta;z),\qquad B_n(\delta;z)=O(z^{2n}).
\end{equation}
If one can take $|\delta|=0$, then we deal with a consistent moment
problem (\ref{SRS}). As we have already mentioned above, this case
is not of big interest for numerical analysis because analytical
methods effectively work. If $|\delta|\ne 0$, then obviously  one
cannot get (\ref{2n-interpolation}) from (\ref{NVZ+}) even for
regular problems. Indeed, given fixed residuals $\delta$,
$|\delta|\ne 0$, the right hand side of (\ref{NVZ+}) is just of order $h_k\delta_k z^k=O(|\delta|)z^k$, $z\to 0$, where $k$ is the
index of the first non-zero residual.

There also exist other serious obstacles to realization of the
 method of small residuals in the interpolation and
approximation problems. The matter is that for inconsistent moment
problems the decreasing of $|\delta|$ to zero always results in
growing to infinity of at least one component of the solution to the
problem (\ref{NVZ}), i.e. the amplitudes $\mu_k(\delta)$ or
frequencies $\lambda_k(\delta)$. This happens because the solution
to the initial, non-regularized inconsistent problem (\ref{SRS})
does not exist. But a similar situation can occur even for some
consistent systems (see Example~\ref{ex6-3}). In such cases the
approximation is impossible: either the computational errors
considerably exceed the residuals or, which is even worse, the
arguments $\lambda_k(\delta) z\to\infty$ leave the domain of the
function $h(z)$ (see Example~\ref{ex6-2}). Moreover, then the
corresponding interpolation formulas are usually in unstable
relation not only to the norms of the residuals but also to their
single components that makes impossible the error estimation via the
norm $|\delta|$ (see Examples~\ref{ex6-1}-\ref{ex6-3}).

We now give several examples of such a divergence for consistent and
inconsistent moment problems. For simplicity, let $n=2$.
\begin{example}
\label{ex6-1} Set $s_0=0$, $s_1=1$, $s_2=0$, $s_3=0$. The system
(\ref{SRS}) is inconsistent. Solving the system (\ref{NVZ}) by the
Prony-Sylvester formulas and Lemma~\ref{400}, we get
$$
\mu_{1}+\mu_2=\delta_0,\qquad \mu_{1}=\frac {(1+\delta_1)^2}{\sqrt
{4\,\delta_3
\delta_1-3\,{\delta_2}^{2}+4\,\delta_3}}+O(|\delta|)\to\infty,
\qquad |\delta|\to 0.
$$
Thus, passage to the limit in $H_n(\delta;z)$ as $|\delta|\to 0$
predictably does not determine any Pad\'{e} amplitude and frequency
sum, moreover, the parameters of the sum obtained are in unstable
relation to the components of the residuals.
\end{example}
\begin{example}
\label{ex6-2} Set $s_0=1$, $s_1=0$, $s_2=0$, $s_3=1$. The system (\ref{SRS})
is inconsistent. By the same formulas for (\ref{NVZ}) we get
 $$
 \lambda_{1}=\frac{1+O\left(|\delta|\right)}{\delta_0\delta_2+\delta_2-\delta_1^2}\;
\to\infty, \qquad |\delta|\to 0,
 $$
i.e. the argument $\lambda_{1} z$ of the amplitude and frequency sum
$H_n(\delta;z)$ in (\ref{NVZ+}) tends to infinity and can leave the
domain of the basis function $h$.
\end{example}
The following example (which arises in the extrapolation problem
considered in Section~\ref{par-extrap}) shows that the method of
small residuals can be unfit even for consistent systems.
\begin{example}
\label{ex6-3} Set $s_0=1$, $s_1=1$, $s_2=1$, $s_3=1$. The system
(\ref{SRS}) is consistent (but not regular); one of its solutions is
obvious: $\lambda_1=\mu_1=1$, $\lambda_2=0$, $\mu_2$ is arbitrary.
However, the method of small residuals leads to the indeterminacy
 $$
 \lambda_1\cdot \lambda_2=\frac {\delta_3+\delta_1+\delta_1\delta_3
 -2\,\delta_2-\delta_2^2}{\delta_2-2\,\delta_1+\delta_0+\delta_0\delta_2-
 \delta_1^{2}}.
$$
It is not clear how to choose the residuals for the convergence of
the process. If one takes, for instance, $\delta_0=\delta_1=0$,
$\delta_2=\delta_3^2$, then again the argument in the amplitude and
frequency sum tends to infinity:
$$
 \lambda_1\cdot \lambda_2= \frac{1-2\delta_3-\delta_3^3}{\delta_3}\;
\to\infty, \qquad |\delta|\to 0.
 $$
\end{example}
\begin{example}
\label{ex6-4} In \cite[\S4.4.2]{YF} the following system arises
in the approximation of the derivative of the Bessel function,
$J'_0$, by amplitude and frequency sums with the basis function
$h(x)=(1-J_0(x))/x$ and $n=2$:
$$
\sum_{k=1}^2\mu_k\lambda_k^{2m+1}=2(m+1),\qquad m=\overline{0,M-1}.
$$
It is easily seen that  it is inconsistent for $M=4$ (cf. the
non-regularized system in the problem of numerical
differentiation in Section~\ref{par-diff}). The authors of \cite{YF}
solve it numerically for $M>4$ and $\varepsilon=2\cdot10^{-16}$ and obtain
the following results:
$$
\mu_1=\overline{\mu_2}=\tfrac{1}{2}+\mu i,\qquad
\lambda_1=\overline{\lambda_2}=1+\lambda i,\qquad \mu\approx
-9.8\cdot 10^7,\qquad \lambda\approx 5.1\cdot 10^{-9}.
$$
It is natural to expect that further decreasing of the residual will
cause the moduli of the amplitudes $\mu_k$ to grow to infinity and
the moduli of the frequencies $\lambda_k$ to tend to one. Again one
gets a divergent process in the context of the Pad\'{e}
interpolation under consideration.
\end{example}

Thus on account of the above-mentioned reasons the method of  small
residuals has to be used for constructing amplitude and frequency
sums with great circumspection as it can lead to unacceptable
results.

Note that the regularization method that we propose also uses
residuals; generally speaking, they are two: $p$ and $q$. The
important difference between our method and the one with small
residuals is in the fact that $p$ and $q$ are fixed, not necessary
small and can be calculated by special formulas depending on the
specificity of the problems considered in
Sections~\ref{regularization_section}-\ref{par-extrap}. For example,
in the problems of numerical differentiation and extrapolation we
obtained explicit expressions for $p$ and $q$, being independent of
$z$, $f$ and $h$. Moreover, the corresponding residuals in the
interpolation disappear not because of decreasing $|\delta|$
but due to adding a fixed regularizing binomial $c_1z^{n-1}+c_2z^{2n-1}$ to an amplitude and frequency
sum. For instance, the following regularized interpolation formula
 corresponds to Example~\ref{ex6-2}:
$$
f(z)=-h_{1}z+\sum_{k=1}^2\mu_kh(\lambda_kz)+O(z^4),\quad
\mu_{1,2}=\tfrac{1}{2}\left(1\mp\tfrac{3\sqrt{5}}{5}\right),\quad
\lambda_{1,2}=-\tfrac{1}{2}\left(1\pm\sqrt{5}\right),
$$
where the remainder depends only on $z$, $f$ and $h$.

In conclusion, it is worth mentioning that for increasing the
rate of approximation one can use the regularization method even for
some regular systems, for example, for those which can be obtained from
non-regular ones by a small variation of the moments.

\section*{Acknowledgments}

We would like to thank the referees for their useful suggestions, which helped to improve the paper.

\bigskip\bigskip\bigskip

\textbf{Petr Chunaev}

Departament de Matem\`{a}tiques, Universitat Aut\`{o}noma de Barcelona

Edifici C, Facultat de Ci\`{e}ncies, 08193 Bellaterra (Barcelona), Spain

e-mail: \textsf{chunayev@mail.ru; chunaev@mat.uab.cat}

\bigskip

\textbf{Vladimir Danchenko}

Functional Analysis and Its Applications Department, Vladimir State University

Belokonskoy str. 3/7, Building 3, 600000 Vladimir, Russia

e-mail: \textsf{vdanch2012@yandex.ru; danch@vlsu.ru}

\end{document}